\documentclass[12pt, reqno, a4paper]{amsart}
\usepackage{amssymb, amsthm, amscd}
\numberwithin{equation}{section}
\newtheorem{state}{Statement}%[theor]
\newtheorem{rem}{Remark}
\newtheorem{example}{Example}
\DeclareMathOperator{\ch}{cosh}
\DeclareMathOperator{\sh}{sinh}
\title[Group Classification of Potential Structures]{Classification
of Potential Structures on Minkowski Space over Subgroups
of the Poincar\'{e} Group}
\author{M.~A.~Parinov}
\address{Ivanovo State University. 153025, Russia, Ivanovo, ul.~Ermaka, 39}
\date{25 February, 2006}
\email{parinov@ivanovo.ac.ru}

\begin{document}
\begin{abstract}
We describe classes of potential structures (covector fields) on Minkowski
space that admit subgroups of the Poincar\'{e} group. We describe also seven
classes of Maxwell spaces that admit subgroups of the Poincar\'{e} group.
\end{abstract}
\maketitle
\medskip
\tableofcontents
%%%%%%%%%%%%%%%%%%%%%%%%%%%%%%%%%%%%%%%%%%%%%%%%%%%%%%%%%%%%%%%%%%%%%%%%
\section{Introduction} % Section 1

Using the classification of subgroups of the Poincar\'{e} group
\cite{Bel'ko} we classified Maxwell spaces with respect these subgroups
\cite{Par03, Par04}. We use classes of potential structures on Minkowski space,
that admit the same subgroups of the Poincar\'{e} group, for obtaining
representatives of Maxwell spaces classes in \cite{Par04}.
Some classes of potential structures were described
in \cite{Polezh-Par, Vorob}, in this paper we present
for the first time the classification of potential structures completely.
This classification is interesting itself, moreover it helps to define more
precisely some classes of Maxwell spaces. For example some classes $C_{p,q}$
in spite of \cite{Par03, Par04} turn out non-empty, we describe them in
appendix.
%%%%%%%%%%%%%%%%%%%%%%%%%%%%%%%%%%%%%%%%%%%%%%%%%%%%%%%%%%%%%%%%%%%%%%%%

\section{Formulation of the problem and method of its solution} % Section 2
\renewcommand{\theequation}{2.\arabic{equation}}

For any smooth, real manifold $M$ we define a {\it potential structure} as a
differential 1-form $A=A_i\,dx^i$, $A_i=A_i(x)$, $x\in M$ \cite{Par03}.

Let $M$ be a four-dimensional manifold; let also $g=g_{ij}dx^idx^j$ be
a pseudo-Euclidean metric on $M$ of Lorentz signature $(---+)$. We
may understand a pair $(M,\,g)$ as a domain in the Minkowski space
$\mathbb R^4_1$. Any triple $(M,\,g,\,A)$ is interpreted as a four-potential
of electromagnetic field. A {\it Maxwell space} is a triple $(M,\,g,\,F)$,
where
\begin{equation}\label{F_{ij}=}
F=dA=F_{ij}dx^idx^j \ \ (F_{ij}=\partial_iA_j-\partial_jA_i)
\end{equation} %\eqno(2.1)
is a generalized symplectic structure \cite{Par03}. Since 2-form $F$ is closed,
$$
dF=0\ \Leftrightarrow\ \partial_{[i} F_{jk]}=0,
$$
then we may understand $F_{ij}$ as a tensor of electromagnetic field%
\footnotemark[1]%
\footnotetext[1]{If the second Maxwell equation
$\nabla_kF^{ik}=-\frac {4\pi}cJ^i$ is satisfied (if we disregard by physical
restrictions, then we may understand this equation as a definition
of current).}.

The problem of group classification of potential structures $(M,\,g,\,A)$
(potentials on $M\subset\mathbb R^4_1$) is analogous to the problem of
classification Maxwell spaces over subgroups of the Poincar\'{e} group
\cite{Par03, Par04}. For every subgroup $G_{p,q}$, corresponding to
the algebra $\mathcal{L}_{p,q}=L\{\xi_1,\,\dots,\,\xi_p\}$\footnotemark[2],%
\footnotetext[2]{See the list of subgroups in \cite{Bel'ko}.}
we find the class $P_{p,q}$ of potentials ${A_i}$, which are invariant
respectively this group; the potential ${A\in P_{p,q}}$ satisfies to
the invariance condition
\begin{equation}\label{CondInv A_i}
L_{\xi_{\alpha}}A_i=0\ \ (\alpha=1,\,\dots,\,p=\text{dim}\mathcal{L}_{p,q})
\end{equation} %\eqno(2.2)
($L_{\xi}$ is the Lie derivative). Solving \eqref{CondInv A_i} for every
algebra $\mathcal{L}_{p,q}$ in \cite{Bel'ko}, we'll get the complete group
classification\footnotemark[3]%
\footnotetext[3]{We execute this operation only for algebras
$\mathcal{L}_{p,q}$ such that $p\le 6$.}
of potential structures.

We take the basis of the Lie algebra corresponding to the Poincar\'{e} group
as follows
$$\aligned
&e_1=(1, 0, 0, 0),\ e_2=(0, 1, 0, 0),\ e_3=(0, 0, 1, 0),\ e_4=(0, 0, 0, 1),\\
&e_{12}=(-x^2, x^1, 0, 0),\ \ e_{13}=(x^3, 0, -x^1, 0),\ \
e_{23}=(0, -x^3, x^2, 0),\\
&e_{14}=(x^4, 0, 0, x^1),\ \ \ e_{24}=(0, x^4, 0, x^2),\ \ \
e_{34}=(0, 0, x^4, x^3).
\endaligned
$$
Here $\{x^i\}$ are the Galilean coordinates such that
$$
g_{ij}=\mathrm{diag}(-1,-1,-1,1).
$$
In what follows, $L\{\xi_1,\,\dots,\,\xi_p\}$ is the linear combination
of vectors $\xi_1,\,\dots,\,\xi_p$. We suppose that components of all tensors
correspond to the Galilean coordinates $\{x^i\}$ even if they are expressed
as functions of other variables.

\begin{rem}
Every potential $A\in P_{p,q}$ admits the group $G_{p,q}$ or more wide
subgroup of the Poincar\'{e} group.
\end{rem}
%%%%%%%%%%%%%%%%%%%%%%%%%%%%%%%%%%%%%%%%%%%%%%%%%%%%%%%%%%%%%%%%%%%%%%%%

\section{Classes of potential structures} % Section 3
\renewcommand{\theequation}{3.\arabic{equation}}
\setcounter{equation}{0}
\subsection{Potentials that admit one-dimensional symmetry groups}
%Subsection  3.1
%
\subsubsection{Translations} % Subsubsection 3.1.1
There are three types of non-conjugate in pairs, one-dimensional subgroups
of translations.

3.1.1.1. {\it Class $P_{1,1a}$.} %% Point 3.1.1.1
The algebra $\mathcal{L}_{1,1a}=L\{e_1\}$ corresponds to the one-dimensional
group $G_{1,1a}$ of translations along the space-like vector $e_1$.
The equation \eqref{CondInv A_i} for the vector $\xi=e_1$ takes the form
\begin{equation}\label{CondInv A_i:e_1}
\partial_1 A_i=0.
\end{equation}%\eqno(3.1)
Therefore all components of covector field $A_i$ are independent of $x^1$.
\begin{state} \label{Cl:P{1,1a}}
The class $P_{1,1a}$ of potentials that admit the group $G_{1,1a}$
consists of the fields $A_i=A_i(x^2,\,x^3,\,x^4)$.
\end{state}

3.1.1.2. {\it Class $P_{1,1b}$.} %% Point 3.1.1.2
The algebra $\mathcal{L}_{1,1b}=L\{e_4\}$ corresponds to the one-dimensional
group $G_{1,1b}$ of translations along the time-like vector $e_4$.
The equation \eqref{CondInv A_i} for the vector $\xi=e_4$ takes the form
\begin{equation}\label{CondInv A_i:e_4}
\partial_4 A_i=0.
\end{equation}%\eqno(3.2)
Therefore all components of covector field $A_i$ are independent of $x^4$.
\begin{state} \label{Cl:P{1,1b}}
The class $P_{1,1b}$ of potentials that admit the group $G_{1,1b}$
consists of the fields $A_i=A_i(x^1,\,x^2,\,x^3)$.
\end{state}

3.1.1.3. {\it Class $P_{1,1c}$.} %% Point 3.1.1.3
The algebra $\mathcal{L}_{1,1c}=L\{e_2+e_4\}$ corresponds to the
one-dimensional group $G_{1,1c}$ of translations along the isotropic vector
$e_2+e_4$. The equation \eqref{CondInv A_i} for the vector $\xi=e_2+e_4$
takes the form
\begin{equation}\label{CondInv A_i:e_2+e_4}
\partial_2 A_i+\partial_4 A_i=0.
\end{equation} %\eqno(3.3)
Using the substitution
\begin{equation}\label{Change:Isotr}
v^1=x^1,\ \ v^2= x^2+x^4,\ \ v^3=x^3,\ \ v^4= x^2- x^4,
\end{equation} %\eqno(3.4)
we receive the solution of equation \eqref{CondInv A_i:e_2+e_4}
\begin{equation}\label{A_i:Isotr}
A_i=A_i(v^1, v^3, v^4)=A_i(x^1, x^3, x^2-x^4),
\end{equation} %\eqno(3.5)
where $A_i(v^1, v^3, v^4)$ are arbitrary functions.
\begin{state} \label{Cl:P{1,1c}}
The class $P_{1,1c}$ of potentials that admit the group $G_{1,1c}$
consists of the fields \eqref{A_i:Isotr}.
\end{state}

\subsubsection{Elliptic helices}  % Subsubsection 3.1.2
The algebra $\mathcal{L}_{1,2}=L\{e_{13}+\lambda e_2+\mu e_4\}$ corresponds to
the one-dimensional group $G_{1,2}$ of elliptic helices or rotations. The
equation \eqref{CondInv A_i} for the vector ${\xi=e_{13}+\lambda e_2+\mu e_4}$
takes the form
\begin{equation}\label{CondInv A_i-el}
x^3\partial_1A_i+\lambda\partial_2A_i-x^1\partial_3A_i+\mu\partial_4A_i+
A_1\delta^3_i-A_3\delta^1_i=0.
\end{equation} %\eqno(3.6)
Using the substitution
$\{x^i\}\mapsto\{\tilde x^i\}=\{r, \tilde x^2, \varphi, \tilde x^4\}$,
\begin{equation}\label{Change:Elipt}
x^1=r\sin\varphi,\ \ x^2=\lambda\varphi+\tilde x^2,\ \
x^3=r\cos\varphi,\ \ x^4=\mu\varphi+\tilde x^4,
\end{equation} %\eqno(3.7)
we transform the equation \eqref{CondInv A_i-el} to the system of equations
\begin{equation}\label{CondInv A_i-el-ch}
\frac{\partial A_1}{\partial\varphi}-A_3=0,\ \
\frac{\partial A_2}{\partial\varphi}=0,\ \
\frac{\partial A_3}{\partial\varphi}+A_1=0,\ \
\frac{\partial A_4}{\partial\varphi}=0.
\end{equation} %\eqno(3.8)
We have the following expression for the solution of the system
\eqref{CondInv A_i-el-ch}
\begin{equation}\label{Cl:Elipt}
\begin{split}
&A_1=C_1\cos\varphi+C_2\sin\varphi,\ \ A_2=A_2(r,\tilde{x}^2,\tilde{x}^4),\\
&A_3=-C_1\sin\varphi+C_2\cos\varphi,\ \ A_4=A_4(r,\tilde{x}^2,\tilde{x}^4),
\end{split}
\end{equation} %\eqno(3.9)
where $C_i=C_i(r,\tilde x^2,\tilde x^4)$ are arbitrary functions.
\begin{state} \label{Cl:P{1,2}}
The class $P_{1,2}$ of potentials that admit the group $G_{1,2}$
consists of the fields \eqref{Cl:Elipt}.
\end{state}

\subsubsection{Hyperbolic helices}  % Subsubsection 3.1.3
The algebra $\mathcal{L}_{1,3}=L\{e_{24}+\lambda e_1\}$ corresponds to the
one-dimensional group $G_{1,3}$ of hyperbolic helices or pseudo-rotations
(Lorentz transformations). The equation \eqref{CondInv A_i} for the vector
$\xi=e_{24}+\lambda e_1$ takes the form
\begin{equation}\label{CondInv A_i-Hyp}
\lambda \partial_1A_i+x^4 \partial_2A_i+x^2 \partial_4A_i+
A_2\delta^4_i+A_4\delta^2_i=0.
\end{equation} %\eqno(3.10)
Using the substitution
\begin{equation}\label{Change:Hyperb1,3-gl2}
x^1=\lambda\varphi+\tilde x^1,\ \ x^2=r\ch\varphi,\ \
x^3=\tilde x^3,\ \ x^4=r\sh\varphi
\end{equation} %\eqno(3.11)
we transform the equation \eqref{CondInv A_i-Hyp} to the system of equations
\begin{equation}\label{CondInv A_i-Hyp-ch}
\frac{\partial A_1}{\partial\varphi}=0,\ \
\frac{\partial A_2}{\partial\varphi}+A_4=0,\ \
\frac{\partial A_3}{\partial\varphi}=0,\ \
\frac{\partial A_4}{\partial\varphi}+A_2=0.
\end{equation} %\eqno(3.12)
We have the following expression for the solution of the system
\eqref{CondInv A_i-Hyp-ch}
\begin{equation}\label{Cl:P{1,3}-gl2}
\begin{split}
&A_1=A_1(\tilde{x}^1,r,\tilde{x}^3),\ \ A_2=C_1\ch\varphi+C_2\sh\varphi,\\
&A_3=A_3(\tilde{x}^1,r,\tilde{x}^3),\ \ A_4=-C_1\sh\varphi-C_2\ch\varphi,
\end{split}
\end{equation} %\eqno(3.13)
where $C_i=C_i(\tilde{x}^1,r,\tilde{x}^3)$ are arbitrary functions.
\begin{state} \label{theor:P{1,3}}
The class $P_{1,3}$ of potentials that admit the group $G_{1,3}$
consists of the fields \eqref{Cl:P{1,3}-gl2}.
\end{state}

\subsubsection{Parabolic helices}  % Subsubsection 3.1.4
The algebra
$\mathcal{L}_{1,4}=L\{e_{12}-e_{14}+\lambda e_2+\mu e_3\}$
($\lambda,\mu=\text{const},\ \lambda\mu=0$)
corresponds to the one-dimensional group $G_{1,4}$ of parabolic helices
or parabolic rotations.
The equation \eqref{CondInv A_i} for the vector
$\xi=e_{12}-e_{14}+\lambda e_2+\mu e_3$ takes the form
\begin{equation}\label{CondInv A_i:prb}
XA_i-A_1(\delta^2_i+\delta^4_i)+(A_2-A_4)\delta^1_i=0,
\end{equation} %\eqno(3.14)
where
\begin{equation}\label{X:prb-gl2}
X=-(x^2+x^4)\partial_1+(x^1+\lambda)\partial_2+\mu\partial_3-x^1\partial_4.
\end{equation} %\eqno(3.15)
We consider 3 cases: a) $\lambda=\mu=0$; b) $\lambda=0$, $\mu\ne0$;
c) $\lambda\ne0$, $\mu=0$.

3.1.4.1. {\it Class $P_{1,4a}$}. %% Point 3.1.4.1
For $\lambda=\mu=0$ we use the substitution
\begin{equation}\label{Change:P{1,4a}}
\begin{split}
&\tilde x^1=x^2+x^4,\ \ \tilde x^2=-{x^1}/(x^2+x^4),\\
&\tilde x^3=x^3,\phantom{i===} \tilde x^4=\frac12(x^1)^2+x^2(x^2+x^4);
\end{split}
\end{equation} %\eqno(3.16)
the operator \eqref{X:prb-gl2} is replaced by partial derivative with respect
to $\tilde x^2$ and the equation \eqref{CondInv A_i:prb} is transformed
to the system of equations
\begin{equation}\label{CondInv A_i:1,4a}
\frac{\partial A_1}{\partial\tilde x^2}+A_2-A_4=0,\,
\frac{\partial A_2}{\partial\tilde x^2}-A_1=0,\,
\frac{\partial A_3}{\partial\tilde x^2}=0,\,
\frac{\partial A_4}{\partial\tilde x^2}-A_1=0.
\end{equation} %\eqno(3.17)
We have the following expression for the solution of the system
\eqref{CondInv A_i:1,4a}
\begin{equation}\label{Cl:P{1,4a}}
\begin{split}
&A_1=C_2\tilde x^2+C_3,\phantom{===} A_2=\frac12C_2(\tilde x^2)^2+
 C_3\tilde x^2+C_1,\\
&A_3=A_3(\tilde x^1,\,\tilde x^3,\,\tilde x^4),\ \ A_4=A_2+C_2,
\end{split}
\end{equation} %\eqno(3.18)
where $A_3(\tilde x^1,\,\tilde x^3,\,\tilde x^4)$ and
$C_k=C_k(\tilde x^1,\,\tilde x^3,\,\tilde x^4)$ $(k=1,\,2,\,3)$
are arbitrary functions.
\begin{state} \label{theor:P{1,4a}}
The class $P_{1,4a}$ of potentials that admit the group $G_{1,4a}$,
corresponding to the algebra $\mathcal{L}_{1,4}$ (${\lambda=\mu=0}$),
consists of the fields \eqref{Cl:P{1,4a}}.
\end{state}

3.1.4.2. {\it Class $P_{1,4b}$}. %% Point 3.1.4.2
For $\lambda=0,\,\mu\ne0$ we use in place of \eqref{Change:P{1,4a}}
the substitution
\begin{equation}\label{Change:P{1,4b}}
\begin{split}
&\tilde x^1=x^2+x^4,\phantom{====} \tilde x^2=-{x^1}/(x^2+x^4),\\
&\tilde x^3=x^3+\frac{\mu x^1}{x^2+x^4},\ \
\tilde x^4=\frac12(x^1)^2+x^2(x^2+x^4);
\end{split}
\end{equation} %\eqno(3.19)
then the equation \eqref{CondInv A_i:prb} is transformed to the system
\eqref{CondInv A_i:1,4a}.
\begin{state} \label{theor:P{1,4b}}
The class $P_{1,4b}$ of potentials that admit the group $G_{1,4b}$,
corresponding to the algebra $\mathcal{L}_{1,4}$ (${\lambda=0}$, ${\mu\ne0}$),
is defined by \eqref{Cl:P{1,4a}}, where the substitution
\eqref{Change:P{1,4a}} is replaced by \eqref{Change:P{1,4b}}.
\end{state}

3.1.4.3. {\it Class $P_{1,4c}$}. %% Point 3.1.4.3
For $\lambda\ne0,\,\mu=0$ we use in place of \eqref{Change:P{1,4a}}
the substitution
\begin{equation}\label{Change:P{1,4c}}
\begin{split}
&\tilde x^1=2\lambda x^1+\left(x^2+x^4\right)^2,\ \
 \tilde x^2=(x^2+x^4)/{\lambda},\ \ \tilde x^3=x^3,\\
&\tilde x^4=\lambda x^4+x^1(x^2+x^4)+\left(x^2+x^4\right)^3/{3\lambda},
\end{split}
\end{equation} %\eqno(3.20)
which transforms \eqref{CondInv A_i:prb} to \eqref{CondInv A_i:1,4a}.
\begin{state} \label{theor:P{1,4c}}
The class $P_{1,4c}$ of potentials that admit the group $G_{1,4c}$,
corresponding to the algebra $\mathcal{L}_{1,4}$ (${\lambda\ne0}$, ${\mu=0}$),
is defined by \eqref{Cl:P{1,4a}}, where the substitution
\eqref{Change:P{1,4a}} is replaced by \eqref{Change:P{1,4c}}.
\end{state}

\subsubsection{Proportional bi-rotations}  % Subsubsection 3.1.5
The algebra $\mathcal{L}_{1,5}=L\{e_{13}+\lambda e_{24}\}$ corresponds
to the group $G_{1,5}$ of proportional bi-rotations. The equation
\eqref{CondInv A_i} for the vector $\xi=e_{13}+\lambda e_{24}$ takes the form
\begin{align}\label{CondInv A_i:birot}
&XA_i+A_1\delta^3_i+\lambda A_2\delta^4_i-A_3\delta^1_i+
 \lambda A_4\delta^2_i,\\ %\eqno(3.21)
             \label{X:birot-gl2}
&X=x^3\partial_1+\lambda x^4\partial_2-x^1\partial_3+
 \lambda x^2\partial_4. %\eqno(3.22)
\end{align}
We use the substitution $\{x^i\}\,\to\,\{r,\,\rho,\,\theta,\,\varphi\}$,
\begin{equation}\label{Change:P{1,5}}
\begin{split}
&x^1=r\cos(\theta-\varphi),\ \ x^2=\rho\ch(\lambda\varphi),\\
&x^3=r\sin(\theta-\varphi),\ \ x^4=\rho\sh(\lambda\varphi);
\end{split}
\end{equation} %\eqno(3.23)
the operator \eqref{X:birot-gl2} is replaced by partial derivative with
respect to $\varphi$ and the equation \eqref{CondInv A_i:birot} is transformed
to the system of equations
\begin{equation}\label{CondInv A_i-birot-ch}
\begin{split}
&\frac{\partial A_1}{\partial\varphi}-A_3=0,\ \
 \frac{\partial A_2}{\partial\varphi}+\lambda A_4=0,\\
&\frac{\partial A_3}{\partial\varphi}+A_1=0,\ \
 \frac{\partial A_4}{\partial\varphi}+\lambda A_2=0.
\end{split}
\end{equation} %\eqno(3.24)
We have the following expression for the solution of the system
\eqref{CondInv A_i-birot-ch}
\begin{equation}\label{Cl:P{1,5}}
\begin{split}
&A_1=C_1\cos\varphi+C_2\sin\varphi,\ \
 A_2=C_3\ch\lambda\varphi+C_4\sh\lambda\varphi,\\
&A_3=-C_1\sin\varphi+C_2\cos\varphi,\ \
 A_4=-C_3\sh\lambda\varphi-C_4\ch\lambda\varphi,
\end{split}
\end{equation} %\eqno(3.25)
where $C_i=C_i(\rho,r,\theta)$ are arbitrary functions.
\begin{state} \label{theor:P_{1,5}}
The class $P_{1,5}$ of potentials that admit the group $G_{1,5}$
consists of the fields \eqref{Cl:P{1,5}}.
\end{state}
%%%-------------------------------------------------------------

\subsection{Potentials that admit two-dimensional symmetry groups}
%Subsection  3.2
%
\subsubsection{Translations}  % Subsubsection 3.2.1
There are three types of non-conjugate in pairs, two-dimensional subgroups
of translations.

3.2.1.1. {\it Class $P_{2,1a}$.} %% Point 3.2.1.1
The algebra $\mathcal{L}_{2,1a}=L\{e_1,\ e_2\}$ corresponds to the group
$G_{2,1a}$ of translations along the vectors of the Euclidean plane.
We have $\mathcal{L}_{1,1a}\subset\mathcal{L}_{2,1a}$, therefore the class
$C_{2,1a}$ is a subclass of the class $C_{1,1a}$.
The equation \eqref{CondInv A_i} for the vector $\xi=e_2$ takes the form
\begin{equation}\label{CondInv A_i:e_2}
\partial_2 A_i=0,
\end{equation}%\eqno(3.26)
Substituting $A_i(x^2,\,x^3,\,x^4)$ for $A_i$ in \eqref{CondInv A_i:e_2},
we get the following result.
\begin{state} \label{state:P{2,1a}}
The class $P_{2,1a}$ of potentials that admit the group $G_{2,1a}$
consists of the fields $A_i=A_i(x^3,\,x^4)$.
\end{state}

3.2.1.2. {\it Class $P_{2,1b}$.} %% Point 3.2.1.2
The algebra $\mathcal{L}_{2,1b}=L\{e_2,\ e_4\}$ corresponds to the group
$G_{2,1b}$ of translations along the vectors of the pseudo-Euclidean plane.
Since $\mathcal{L}_{1,1b}\subset\mathcal{L}_{2,1b}$, we have
$P_{2,1b}\subset P_{1,1b}$.
Substituting $A_i(x^1,\,x^2,\,x^3)$ for $A_i$ in \eqref{CondInv A_i:e_2},
we get the following result.
\begin{state} \label{state:P{2,1b}}
The class $P_{2,1b}$ of potentials that admit the group $G_{2,1b}$
consists of the fields $A_i=A_i(x^1,\,x^3)$.
\end{state}

3.2.1.3. {\it Class $P_{2,1c}$.} %% Point 3.2.1.3
The algebra $\mathcal{L}_{2,1c}=L\{e_1,\ e_2+e_4\}$ corresponds to the group
$G_{2,1c}$ of translations along the vectors of the isotropic plane.
Since $\mathcal{L}_{1,1c}\subset\mathcal{L}_{2,1c}$, we have
$P_{2,1c}\subset P_{1,1c}$.
Combining \eqref{A_i:Isotr} and \eqref{CondInv A_i:e_1}, we get
the following result.
\begin{state} \label{state:P{2,1c}}
The class $P_{2,1c}$ of potentials that admit the group $G_{2,1c}$
consists of the fields $A_i=A_i(x^3,\,x^2-x^4)$.
\end{state}

\subsubsection{Class $P_{2,2}$}  % Subsubsection 3.2.2
The algebra $\mathcal{L}_{2,2}=L\{e_{13}+\mu e_4,\ e_2\}$ corresponds to
the group $G_{2,2}$ generated by elliptic helices with a time-like axis
and by translations along a space-like straight line. Let
$\mathcal{L}_{1,2b}=L\{e_{13}+\mu e_4\}$ be the algebra $\mathcal{L}_{1,2}$
for $\lambda=0$. Since $\mathcal{L}_{1,2b}\subset\mathcal{L}_{2,2}$, we have
$P_{2,2}\subset P_{1,2b}$. Since $\lambda=0$, the substitution
\eqref{Change:Elipt} takes the form
\begin{equation}\label{Change:Elipt-b}
x^1=r\sin\varphi,\ \ x^2=\tilde x^2,\ \
x^3=r\cos\varphi,\ \ x^4=\mu\varphi+\tilde x^4.
\end{equation} %\eqno(3.27)
Substituting \eqref{Cl:Elipt} for $A_i$ in \eqref{CondInv A_i:e_2}, we get
\begin{equation}\label{Cl:P{2,2}}
\begin{split}
&A_1=C_1\cos\varphi+C_2\sin\varphi,\ \ A_2=A_2(r,\tilde{x}^4),\\
&A_3=-C_1\sin\varphi+C_2\cos\varphi,\ \ A_4=A_4(r,\tilde{x}^4),
\end{split}
\end{equation} %\eqno(3.28)
where $C_i=C_i(r,\tilde x^4)$ are arbitrary functions.
\begin{state} \label{state:P{2,2}}
The class $P_{2,2}$ of potentials that admit the group $G_{2,2}$
consists of the fields \eqref{Cl:P{2,2}}.
\end{state}

\subsubsection{Class $P_{2,3}$}  % Subsubsection 3.2.3
The algebra $\mathcal{L}_{2,3}=L\{e_{13}+\lambda e_2, e_4\}$ corresponds to
the group $G_{2,3}$ generated by elliptic helices with a space-like axis
and by translations along a time-like straight line.
By $P_{1,2a}$ denote the class $P_{1,2}$ for $\mu=0,\,\lambda\ne0$.
The class $P_{2,3}$ is an intersection of $P_{1,1b}$ and $P_{1,2a}$.
In this case, the substitution \eqref{Change:Elipt} takes the form
\begin{equation} \label{Change:Elipt-a}
x^1=r\sin\varphi,\ \ x^2=\lambda\varphi+\tilde x^2,\ \
x^3=r\cos\varphi,\ \ x^4=\tilde x^4.
\end{equation} %\eqno(3.29)
Since $P_{2,3}\subset P_{1,2}$, it follows that $P_{2,3}$ is defined by
\eqref{Cl:Elipt}. If $P_{2,3}\subset P_{1,1b}$, then all components $A_i$
are independent of $x^4$, therefore we have the following result.
\begin{state} \label{state:P{2,3}}
The class $P_{2,3}$ of potentials that admit the group $G_{2,3}$
consists of the following fields
\begin{equation}\label{Cl:P{2,3}}
\begin{split}
&A_1=b_1\cos\varphi+b_2\sin\varphi,\ \ A_2=A_2(r,\tilde x^2),\\
&A_3=-b_1\sin\varphi+b_2\cos\varphi,\ \ A_4=A_4(r,\tilde x^2),
\end{split}
\end{equation} %\eqno(3.30)
where $b_k=b_k(r,\tilde x^2)$, $A_2(r,\tilde x^2)$, and $A_4(r,\tilde x^2)$
are arbitrary functions (the transformation of coordinates is defined by
\eqref{Change:Elipt-a}).
\end{state}

\subsubsection{Class $P_{2,4}$} %% Subsubsection 3.2.4
The algebra $\mathcal{L}_{2,4}=L\{e_{13}+\lambda e_2,\ e_2+e_4\}$ corresponds
to the group $G_{2,4}$ generated by elliptic helices with a space-like axis
and by translations along an isotropic straight line.
The class $P_{2,4}$ is an intersection of classes $P_{1,2a}$ and $P_{1,1c}$.
Substituting \eqref{Cl:Elipt} for $A_i$ in \eqref{CondInv A_i:e_2+e_4}, we get
the following result.
\begin{state} \label{state:P{2,4}}
The class $P_{2,4}$ of potentials that admit the group $G_{2,4}$
consists of the following fields
\begin{equation}\label{Cl:P{2,4}}
\begin{split}
&A_1=C_1\cos\varphi+C_2\sin\varphi,\ \ A_2=A_2(r,\tilde{x}^2-\tilde{x}^4),\\
&A_3=-C_1\sin\varphi+C_2\cos\varphi,\ \ A_4=A_4(r,\tilde{x}^2-\tilde{x}^4),
\end{split}
\end{equation} %\eqno(3.31)
where $C_i=C_i(r,\tilde x^2-\tilde x^4)$ are arbitrary functions (the
transformation of coordinates is defined by \eqref{Change:Elipt-a}).
\end{state}

\subsubsection{Class $P_{2,5}$} %% Subsubsection 3.2.5
The algebra $\mathcal{L}_{2,5}=L\{e_{24}+\lambda e_3,\,e_1\}$ corresponds to
the group $G_{2,5}$ generated by hyperbolic helices and by translations
along the space-like straight line.
By $P_{1,3a}$ denote the class of potentials that admit the group $G_{1,3a}$
corresponding to the algebra
$$
\mathcal{L}_{1,3a}=L\{e_{24}+\lambda e_3\};
$$
it is defined by \eqref{Cl:P{1,3}-gl2}, where the substitution
\eqref{Change:Hyperb1,3-gl2} is replaced by the following one:
\begin{equation}\label{Change:Hyperb1,3a-gl2}
x^1=\tilde x^1,\ \ x^2=r\ch\varphi,\ \
x^3=\lambda\varphi+\tilde x^3,\ \ x^4=r\sh\varphi.
\end{equation} %\eqno(3.32)
The class $P_{2,5}$ is an intersection of classes $P_{1,3a}$ and $P_{1,1a}$.
Combining \eqref{Cl:P{1,3}-gl2}, \eqref{Change:Hyperb1,3a-gl2}, and
\eqref{CondInv A_i:e_1}, we get the following result.
\begin{state}\label{state:P{2,5}}
The class $P_{2,5}$ of potentials that admit the group $G_{2,5}$
consists of the following fields
\begin{equation}\label{Cl:P{2,5}}
\begin{split}
&A_1=A_1(r,\tilde{x}^3),\ \ A_2=C_1\ch\varphi+C_2\sh\varphi,\\
&A_3=A_3(r,\tilde{x}^3),\ \ A_4=-C_1\sh\varphi-C_2\ch\varphi,
\end{split}
\end{equation} %\eqno(3.33)
where $C_i=C_i(r,\tilde{x}^3)$ are arbitrary functions and the transformation
of coordinates is defined by \eqref{Change:Hyperb1,3a-gl2}.
\end{state}

\subsubsection{Class $P_{2,6}$} %% Subsubsection 3.2.6
The algebra $\mathcal{L}_{2,6}=L\{e_{24}+\lambda e_3,\,e_2-e_4\}$ corresponds
to the group $G_{2,6}$ generated by hyperbolic helices and by translations
along the isotropic straight line.
Since $\mathcal{L}_{1,3a}\subset\mathcal{L}_{2,6}$, we have
${P_{2,6}\subset P_{1,3a}}$. The equation \eqref{CondInv A_i} for the vector
$\xi=e_2-e_4$ takes the form
\begin{equation}\label{CondInv A_i:e_2-e_4}
\partial_2A_i-\partial_4A_i=0.
\end{equation} %\eqno(3.34)
Using the substitution \eqref{Change:Hyperb1,3a-gl2}, we solve the equation
\eqref{CondInv A_i:e_2-e_4} for the potential \eqref{Cl:P{1,3}-gl2}.
We get the following result.
\begin{state}\label{state:P{2,6}}
The class $P_{2,6}$ of potentials that admit the group $G_{2,6}$
consists of the following fields
\begin{equation}\label{Cl:P{2,6}}
\begin{split}
&A_1=A_1(\tilde{x}^1,\tilde{x}^3-\lambda\ln r),\ \
 A_2=C_1\ch\varphi+C_2\sh\varphi,\\
&A_3=A_3(\tilde{x}^1,\tilde{x}^3-\lambda\ln r),\ \
 A_4=-C_1\sh\varphi-C_2\ch\varphi,
\end{split}
\end{equation} %\eqno(3.35)
where
\begin{equation}\label{Cl:P{2,6}-a}
\begin{split}
&C_1=a_1(\tilde{x}^1,\tilde{x}^3-\lambda\ln r)\ch\ln r+
     a_2(\tilde{x}^1,\tilde{x}^3-\lambda\ln r)\sh\ln r,\\
&C_2=a_1(\tilde{x}^1,\tilde{x}^3-\lambda\ln r)\sh\ln r+
     a_2(\tilde{x}^1,\tilde{x}^3-\lambda\ln r)\ch\ln r;
\end{split}
\end{equation} %\eqno(3.36)
the transformation of coordinates is defined by \eqref{Change:Hyperb1,3a-gl2}.
\end{state}

\subsubsection{} %% Subsubsection 3.2.7
Here we describe classes of potentials corresponding to the algebra
$\mathcal{L}_{2,7}=L\{e_{12}-e_{14}+\lambda e_2+\mu e_3,\ e_2-e_4\}\ \
(\lambda\mu=0)$ for various $\lambda$ and $\mu$. The corresponding group
$G_{2,7}$ is generated by parabolic helices and by translations along the
isotropic straight line. The algebra $\mathcal{L}_{2,7}$ is an extension of
$\mathcal{L}_{1,4}$ by means of the vector $\xi=e_2-e_4$; therefore
corresponding classes $P_{2,7a}$, $P_{2,7b}$, and $P_{2,7c}$ are restrictions
of classes $P_{1,4a}$, $P_{1,4b}$, and $P_{1,4c}$ by the condition
\eqref{CondInv A_i:e_2-e_4}.
\smallskip

3.2.7.1. {\it Class $P_{2,7a}$.} %% Point 3.2.7.1
For $\lambda=\mu=0$ the substitution \eqref{Change:P{1,4a}} transforms the
equation \eqref{CondInv A_i:e_2-e_4} to the form
$\tilde x^1\, {\partial A_i}/{\partial\tilde x^4}=0$; therefore all
components of potential are independent of $\tilde x^4$.
\begin{state} \label{theor:P{2,7a}}
The class $P_{2,7a}$ of potentials that admit the group $G_{2,7a}$,
corresponding to the algebra $\mathcal{L}_{2,7}$ (${\lambda=\mu=0}$),
consists of the following fields
\begin{equation}\label{Cl:P{2,7a}}
\begin{split}
&A_1=C_2\tilde x^2+C_3,\ \ A_2=\frac12C_2(\tilde x^2)^2+ C_3\tilde x^2+C_1,\\
&A_3=A_3(\tilde x^1,\,\tilde x^3),\ \ A_4=A_2+C_2,
\end{split}
\end{equation} %\eqno(3.37)
where $A_3(\tilde x^1,\,\tilde x^3)$ and $C_k=C_k(\tilde x^1,\,\tilde x^3)$
$(k=1,\,2,\,3)$ are arbitrary functions and the transformation of coordinates
is defined by \eqref{Change:P{1,4a}}.
\end{state}

3.2.7.2. {\it Class $P_{2,7b}$.} %% Point 3.2.7.2
For $\lambda=0,\,\mu\ne0$ we use the substitution \eqref{Change:P{1,4b}}
instead of \eqref{Change:P{1,4a}}. We have the following result.
\begin{state} \label{theor:P{2,7b}}
The class $P_{2,7b}$ of potentials that admit the group $G_{2,7b}$,
corresponding to the algebra $\mathcal{L}_{2,7}$ ${(\lambda=0}$,\ ${\mu\ne0)}$,
consists of the fields \eqref{Cl:P{2,7a}}, where the transformation
of coordinates is defined by \eqref{Change:P{1,4b}}.
\end{state}

3.2.7.3. {\it Class $P_{2,7c}$.} %% Point 3.2.7.3
For $\lambda\ne0$ and $\mu=0$ the substitution \eqref{Change:P{1,4c}}
transforms the equation \eqref{CondInv A_i:e_2-e_4} to the form:
$-\lambda\,{\partial A_i}/{\partial\tilde x^4}=0$; therefore all
components of potential are independent of $\tilde x^4$.
\begin{state} \label{theor:P{2,7c}}
The class $P_{2,7c}$ of potentials that admit the group $G_{2,7c}$,
corresponding to the algebra $\mathcal{L}_{2,7}$ ${(\lambda\ne0,}$ ${\mu=0)}$,
consists of the fields \eqref{Cl:P{2,7a}}, where the transformation
of coordinates is defined by \eqref{Change:P{1,4c}}.
\end{state}

\subsubsection{Class $P_{2,8}$} %% Subsubsection 3.2.8
The algebra $\mathcal{L}_{2,8}=L\{e_{12}-e_{14}+\lambda e_2,\ e_3\}$
corresponds to the group $G_{2,8}$ generated by parabolic helices and by
translations along a space-like straight line. Since
$\mathcal{L}_{1,4c}\subset\mathcal{L}_{2,8}$, then ${P_{2,8}\subset P_{1,4c}}$.
The class $P_{2,8}$ is a restriction of the class $P_{1,4c}$ by the condition
\eqref{CondInv A_i} for the vector $e_3$
\begin{equation}\label{CondInv A_i:e_3}
\partial_3A_i=0.
\end{equation} %\eqno(3.38)
The substitution \eqref{Change:P{1,4c}} transforms the equation
\eqref{CondInv A_i:e_3} to the form: ${\partial A_i/\partial\tilde x^3=0}$;
thus we have the following result.
\begin{state} \label{theor:P{2,8}}
The class $P_{2,8}$ of potentials that admit the group $G_{2,8}$
consists of the following fields
\begin{equation}\label{Cl:P{2,8}}
\begin{split}
&A_1=C_2\tilde x^2+C_3,\ \ A_2=\frac12C_2(\tilde x^2)^2+
 C_3\tilde x^2+C_1,\\
&A_3=A_3(\tilde x^1,\,\tilde x^4),\ \ A_4=A_2+C_2,
\end{split}
\end{equation} %\eqno(3.39)
where $A_3(\tilde x^1,\,\tilde x^4)$ and $C_k=C_k(\tilde x^1,\,\tilde x^4)$
$(k=1,\,2,\,3)$ are arbitrary functions and the transformation
of coordinates is defined by \eqref{Change:P{1,4c}}.
\end{state}

\subsubsection{Class $P_{2,9}$} %% Subsubsection 3.2.9
The algebra $\mathcal{L}_{2,9}=L\{e_{13}+\lambda e_{24},\,e_2-e_4\}$
corresponds to the group $G_{2,9}$ generated by proportional bi-rotations
and by translations along an isotropic straight line.
Since $\mathcal{L}_{1,5}\subset\mathcal{L}_{2,9}$, then the class $P_{2,9}$
is a subclass of the class $P_{1,5}$.
The class $P_{2,9}$ is a restriction of the class $P_{1,5}$ by the condition
\eqref{CondInv A_i:e_2-e_4} (\eqref{CondInv A_i} for the vector $e_2-e_4$).
The substitution \eqref{Change:P{1,5}} transforms the equation
\eqref{CondInv A_i:e_2-e_4} to the form:
\begin{equation}\label{CondInv A_i:e_2-e_4-Ch1,5}
\frac{\partial A_i}{\partial\varphi}+\frac{\partial A_i}{\partial\theta}-
\lambda\rho\frac{\partial A_i}{\partial\rho}=0;
\end{equation} %\eqno(3.40)
thus we have the following result.
\begin{state} \label{theor:P_{2,9}}
The class $P_{2,9}$ of potentials that admit the group $G_{2,9}$
consists of the following fields
\begin{equation}\label{Cl:P{2,9}-1}
\begin{split}
&A_1=C_1\cos\varphi+C_2\sin\varphi,\ \
 A_2=\rho\Phi_3e^{\lambda\varphi},\\
&A_3=-C_1\sin\varphi+C_2\cos\varphi,\ \
 A_4=-\rho\Phi_3e^{\lambda\varphi},
\end{split}
\end{equation} %\eqno(3.41)
\begin{equation}\label{Cl:P{2,9}-2}
\begin{split}
&C_1=\Phi_1\cos\frac{\ln\rho}{\lambda}+\Phi_2\sin\frac{\ln\rho}{\lambda},\\
&C_2=-\Phi_1\sin\frac{\ln\rho}{\lambda}+\Phi_2\cos\frac{\ln\rho}{\lambda},
\end{split}
\end{equation} %\eqno(3.42)
where $\Phi_k=\Phi_k(r,\,\lambda\theta+\ln\rho)$ are arbitrary functions
and the transformation of coordinates is defined by \eqref{Change:P{1,5}}.
\end{state}

\subsubsection{Class $P_{2,10}$} %% Subsubsection 3.2.10
The algebra
$\mathcal{L}_{2,10}=L\{e_{13},\,e_{24}\}=L\{e_{13}+e_{24},\,e_{13}\}$
corresponds to the group $G_{2,10}$ generated by rotations and pseudo-rotations
or, equivalently, by rotations and proportional bi-rotations for $\lambda=1$.
As $\mathcal{L}_{1,5}\subset\mathcal{L}_{2,10}$, then
$C_{2,10}\subset C_{1,5}$ ($\lambda=1$). In this case \eqref{Cl:P{1,5}} takes
the form
\begin{equation}\label{Cl:P{1,5}'}
\begin{split}
&A_1=C_1\cos\varphi+C_2\sin\varphi,\ \ \
 A_2=C_3e^{\varphi}+C_4e^{-\varphi},\\
&A_3=-C_1\sin\varphi+C_2\cos\varphi,\ \
 A_4=-C_3e^{\varphi}+C_4e^{-\varphi}.
\end{split}
\end{equation} %\eqno(3.43)
Substituting \eqref{Cl:P{1,5}'} for $A_i$ in the equation
\eqref{CondInv A_i-el} for $\lambda=\mu=0$
\begin{equation}\label{CondInv A_i:e_{13}}
x^3\partial_1A_i-x^1\partial_3A_i+A_1\delta^3_i-A_3\delta^1_i=0,
\end{equation} %\eqno(3.44)
we obtain the following result.
\begin{state} \label{theor:P_{2,10}}
The class $P_{2,10}$ of potentials that admit the group $G_{2,10}$
consists of the following fields
\begin{equation}\label{Cl:P{2,10}}
\begin{split}
&A_1=-t_1\sin(\theta-\varphi)+t_2\cos(\theta-\varphi),\ \
 A_2=t_3e^{\varphi}+t_4e^{-\varphi},\\
&A_3=t_1\cos(\theta-\varphi)+t_2\sin(\theta-\varphi),\ \
 A_4=-t_3e^{\varphi}+t_4e^{-\varphi},
\end{split}
\end{equation} %\eqno(3.45)
where $t_k=t_k(r,\,\rho)$ are arbitrary functions and the transformation
of coordinates is defined by \eqref{Change:P{1,5}} for $\lambda=1:$
\begin{equation}\label{Change:P{2,10}}
\begin{split}
&x^1=r\cos(\theta-\varphi),\ \ x^2=\rho\ch\varphi,\\
&x^3=r\sin(\theta-\varphi),\ \ x^4=\rho\sh\varphi.
\end{split}
\end{equation} %\eqno(3.46)
\end{state}

\subsubsection{}  %% Subsubsection 3.2.11
Here we describe classes of potentials corresponding to the algebra
$\mathcal{L}_{2,11}=
L\{e_{12}-e_{14}+\lambda e_1+\mu e_3,\ e_{23}+e_{34}-\mu e_1+\lambda e_3\}$
\linebreak
$({\lambda=0},\,\mu\ne0 \sim \lambda\ne0,\,\mu=0)$. The case $\lambda=\mu=0$
is required for description some following classes.

3.2.11.1. {\it Class $P_{2,11}$ $(\lambda=0,\,\mu\ne0)$.} %% Point 3.2.11.1
The algebra $\mathcal{L}_{2,11}$ corresponds to the group $G_{2,11}$ generated
by two one-dimensional subgroups of parabolic helices with different axises.
Since $\mathcal{L}_{1,4b}\subset\mathcal{L}_{2,11}$, then
$P_{2,11}\subset P_{1,4b}$. For description  the class $P_{2,11}$ we
substitute \eqref{Cl:P{1,4a}}--\eqref{Change:P{1,4b}} for $A_i$ in equation
\eqref{CondInv A_i} for the vector $\xi=e_{23}+e_{34}-\mu e_1$
\begin{align}\label{CondInv A_i:e_{23}+e_{34}-mu e_1}
&XA_i-A_2\delta^3_i+A_3(\delta^2_i+\delta^4_i)+A_4\delta^3_i=0,\\ %\eqno(3.47)
              \label{X:e_{23}+e_{34}-mu e_1:gl2}
Xf&=-\mu\partial_1f-x^3\partial_2f+(x^2+x^4)\partial_3f+
   x^3\partial_4f=\nonumber\\
  &=\frac{\mu}{\tilde x^1}\frac{\partial f}{\partial\tilde x^2}+
\frac{(\tilde x^1)^2-\mu^2}{\tilde x^1}\frac{\partial f}{\partial\tilde x^3}-
\tilde x^1\tilde x^3\frac{\partial f}{\partial\tilde x^4}. %\eqno(3.48)
\end{align}
We use the substitution
\begin{equation}\label{Cl:P_{2,11}-3:gl2}
\begin{split}
u=\frac{\tilde x^1\tilde x^3}{(\tilde x^1)^2-\mu^2},\ \
v=\frac12(\tilde x^1\tilde x^3)^2+\tilde x^4\left((\tilde x^1)^2-\mu^2\right)
\end{split}
\end{equation} %\eqno(3.49)
to solve these equations; we obtain the following result.
\begin{state} \label{theor:P_{2,11}}
The class $P_{2,11}$ of potentials that admit the group $G_{2,11}$
consists of the following fields
\begin{equation}\label{Cl:P_{2,11}-1:gl2}
\begin{split}
&A_1=\Phi\tilde x^2+\Psi,\ \
 A_2=\frac12\Phi\left(\tilde x^2\right)^2+\Psi\tilde x^2+\Xi,\\
&A_3=\Upsilon,\ \
 A_4=\frac12\Phi\left(\tilde x^2\right)^2+\Psi\tilde x^2+\Xi+\Phi,
\end{split}
\end{equation} %\eqno(3.50)
\begin{equation}\label{Cl:P_{2,11}-2:gl2}
\begin{split}
&\Psi=-\frac{\mu u}{\tilde x^1}\Phi+C_1,\ \ \Upsilon=-\Phi u+C_2,\\
&\Xi=\frac{\mu^2+\left(\tilde x^1\right)^2}
{2\left(\tilde x^1\right)^2}\Phi u^2-
\frac{\mu C_1+\tilde x^1C_2}{\tilde x^1}u+C_3,
\end{split}
\end{equation} %\eqno(3.51)
where $\Phi=\Phi(\tilde x^1,\,v)$ and $C_k=C_k(\tilde x^1,\,v)$
are arbitrary functions and transformations of variables are defined by
\eqref{Change:P{1,4b}} and \eqref{Cl:P_{2,11}-3:gl2}.
\end{state}

3.2.11.2. {\it Class $P_{2,11a}$ $(\lambda=\mu=0)$.} %% Point 3.2.11.2
The algebra
$$
\mathcal{L}_{2,11a}=L\{e_{12}-e_{14},\,e_{23}+e_{34}\}
$$
corresponds to the group $G_{2,11a}$ generated by two one-dimensional
subgroups of parabolic rotations. We write the equation \eqref{CondInv A_i}
for the basis vectors $e_{12}-e_{14}$ and $e_{23}+e_{34}$:
\begin{align}\label{CondInv A_i:e_{12}-e_{14}(2,11a)}
&XA_1+A_2-A_4=0,\ \ XA_2-A_1=0,\ \ XA_3=0,\ \ XA_4-A_1=0,\nonumber\\
&X=-(x^2+x^4)\partial_1+x^1\partial_2-x^1\partial_4,
\end{align} %\eqno(3.52)
and
\begin{align}\label{CondInv A_i:e_{23}+e_{34}(2,11a)}
&YA_1=0,\ \ YA_2+A_3=0,\ \ YA_3-A_2+A_4=0,\ \ YA_4+A_3=0,\nonumber\\
&Y=-x^3\partial_2+(x^2+x^4)\partial_3+x^3\partial_4.
\end{align} %\eqno(3.53)
We use the substitution
\begin{equation}\label{Change:P_{2,11a}}
\begin{split}
&\tilde x^1=x^2+x^4,\ \ \tilde x^2=-{x^1}/(x^2+x^4),\ \
 \tilde x^3={x^3}/(x^2+x^4),\\
&\tilde x^4=(x^1)^2+(x^2)^2+(x^3)^2-(x^4)^2;
\end{split}
\end{equation} %\eqno(3.54)
the operator $X$ is replaced by partial derivative with respect to
$\tilde x^2$, the operator $Y$ --- by partial derivative with respect to
$\tilde x^3$. We have the following solution of the system
\eqref{CondInv A_i:e_{12}-e_{14}(2,11a)}--%
\eqref{CondInv A_i:e_{23}+e_{34}(2,11a)}:
\begin{align}\label{Solution:P_{2,11a}}
&A_1=-\tilde x^2\,\Phi+\Psi,\ \
 A_2=-\frac12\Phi((\tilde x^2)^2+(\tilde x^3)^2)+
\tilde x^2\,\Psi-\tilde x^3\Xi+\Theta,\nonumber\\
&A_3=\tilde x^3\,\Phi+\Xi,\ \ \ \ A_4=A_2-\Phi,
\end{align} %\eqno(3.55)
where $\Phi=\Phi(\tilde x^1,\tilde x^4)$, $\Psi=\Psi(\tilde x^1,\tilde x^4)$,
$\Xi=\Xi(\tilde x^1,\tilde x^4)$ and $\Theta=\Theta(\tilde x^1,\tilde x^4)$
are arbitrary functions.
\begin{state} \label{theor:P_{2,11a}}
The class $P_{2,11a}$ of potentials that admit the group $G_{2,11a}$ is
defined by \eqref{Solution:P_{2,11a}} and \eqref{Change:P_{2,11a}}.
\end{state}

\subsubsection{Class $P_{2,12}$}  %% Subsubsection 3.2.12
The algebra $\mathcal{L}_{2,12}=L\{e_{12}-e_{14},\ e_{24}+\lambda e_3\}$
corresponds to the group $G_{2,12}$ generated by parabolic rotations and
by hyperbolic helices.
Since $\mathcal{L}_{1,4a}\subset\mathcal{L}_{2,12}$ and
$\mathcal{L}_{1,3a}\subset\mathcal{L}_{2,12}$, then
$P_{2,12}=P_{1,4a}\cap P_{1,3a}$. For description  the class $P_{2,12}$ we
substitute \eqref{Cl:P{1,4a}}--\eqref{Change:P{1,4a}} for $A_i$ in equation
\eqref{CondInv A_i} for the vector $\xi=e_{24}+\lambda e_3$
\begin{equation}\label{CondInv A_i:e_{24}+lambda e_3}
x^4\partial_2A_i+\lambda\partial_3A_i+x^2\partial_4A_i+
A_2\delta^4_i+A_4\delta^2_i=0.
\end{equation} %\eqno(3.56)
The equation \eqref{CondInv A_i:e_{24}+lambda e_3} is a system
\begin{equation}\label{CondInv A_i:e_{24}+lambda e_3-X}
XA_1=0,\ \ XA_2+A_4=0,\ \ XA_3=0,\ \ XA_4+A_2=0,
\end{equation} %\eqno(3.57)
where the operator $X$ by substitution \eqref{Change:P{1,4a}} is replaced
to the form:
$$
Xf=x^4\partial_2f+\lambda\partial_3f+x^2\partial_4f=
\tilde x^1\frac{\partial f}{\partial\tilde x^1}-
\tilde x^2\frac{\partial f}{\partial\tilde x^2}+
\lambda\frac{\partial f}{\partial\tilde x^3}+
(\tilde x^1)^2\frac{\partial f}{\partial\tilde x^4}.
$$
Substituting \eqref{Cl:P{1,4a}} for $A_i$ in
\eqref{CondInv A_i:e_{24}+lambda e_3-X},
we obtain some differential equation; taking into account a linear
independence of the functions $(\tilde x^2)^2$, $\tilde x^2$, and 1, we get
the following system
\begin{equation}\label{Cl:P_{2,12}-sist}
\begin{split}
&\tilde x^1\frac{\partial C_2}{\partial\tilde x^1}+
 \lambda\frac{\partial C_2}{\partial\tilde x^3}+
 (\tilde x^1)^2\frac{\partial C_2}{\partial\tilde x^4}-C_2=0,\\
&\tilde x^1\frac{\partial C_3}{\partial\tilde x^1}+
 \lambda\frac{\partial C_3}{\partial\tilde x^3}+
 (\tilde x^1)^2\frac{\partial C_3}{\partial\tilde x^4}=0,\\
&\tilde x^1\frac{\partial A_3}{\partial\tilde x^1}+
 \lambda\frac{\partial A_3}{\partial\tilde x^3}+
 (\tilde x^1)^2\frac{\partial A_3}{\partial\tilde x^4}=0,\\
&\tilde x^1\frac{\partial C_1}{\partial\tilde x^1}+
 \lambda\frac{\partial C_1}{\partial\tilde x^3}+
 (\tilde x^1)^2\frac{\partial C_1}{\partial\tilde x^4}+C_1+C_2=0,\\
&\tilde x^1\frac{\partial C_1}{\partial\tilde x^1}+
 \lambda\frac{\partial C_1}{\partial\tilde x^3}+
 (\tilde x^1)^2\frac{\partial C_1}{\partial\tilde x^4}+C_1+\\
&+\tilde x^1\frac{\partial C_2}{\partial\tilde x^1}+
 \lambda\frac{\partial C_2}{\partial\tilde x^3}+
 (\tilde x^1)^2\frac{\partial C_2}{\partial\tilde x^4}=0.
\end{split}
\end{equation} %\eqno(3.58)
Using the substitution
\begin{equation}\label{Change:P_{2,12}}
u=\tilde x^3-\lambda\ln\tilde x^1,\ \
v=\tilde x^4-\frac12\left(\tilde x^1\right)^2,
\end{equation} %\eqno(3.59)
we integrate the system \eqref{Cl:P_{2,12}-sist}; the result is
\begin{equation}\label{Cl:P_{2,12}-func}
C_1=-\frac{\tilde x^1}2\Phi_1+\frac1{\tilde x^1}\Phi_2,\ \
C_2=\tilde x^1\Phi_1,\ \ C_3=\Phi_3,\ \ A_3=\Phi_4,
\end{equation} %\eqno(3.60)
where $\Phi_k=\Phi_k(u,\,v)$ are arbitrary functions.
Substituting \eqref{Cl:P_{2,12}-func} for $C_k$ and $A_3$ in
\eqref{Cl:P{1,4a}}, we obtain the following result.
\begin{state} \label{theor:P_{2,12}}
The class $P_{2,12}$ of potentials that admit the group $G_{2,12}$ is
defined by
\begin{equation}\label{Cl:P_{2,12}}
\begin{split}
&A_1=\tilde x^1\tilde x^2\Phi_1+\Phi_3,\,
 A_2=\frac{\tilde x^1}2\left(\left(\tilde x^2\right)^2-1\right)\Phi_1+
     \tilde x^2\Phi_3+\frac{\Phi_2}{\tilde x^1},\\
&A_3=\Phi_4,\phantom{=====i}
 A_4=\frac{\tilde x^1}2\left(\left(\tilde x^2\right)^2+1\right)\Phi_1+
     \tilde x^2\Phi_3+\frac{\Phi_2}{\tilde x^1},
\end{split}
\end{equation} %\eqno(3.61)
where transformations of variables are defined by
\eqref{Change:P{1,4a}} and \eqref{Change:P_{2,12}}.
\end{state}
%%%-------------------------------------------------------------

\subsection{Potentials that admit three-dimensional symmetry groups}
%Subsection  3.3
%
\subsubsection{}  % Subsubsection 3.3.1
Here we describe classes of potentials corresponding to three-dimensional
groups of translations.

3.3.1.1. {\it Class $P_{3,1a}$} %% Point 3.3.1.1
The algebra $\mathcal{L}_{3,1a}=L\{e_1,e_2,e_3\}$ corresponds to the group
$G_{3,1a}$ of translations along the vectors of three-di\-men\-si\-onal
Euclidean space $Ox^1x^2x^3$. Since
$\mathcal{L}_{2,1a}\subset\mathcal{L}_{3,1a}$, then the class $P_{3,1a}$
is a subclass of $P_{2,1a}$. Substituting $A_i(x^3,\,x^4)$ for $A_i$ in
the equation \eqref{CondInv A_i:e_3} (\eqref{CondInv A_i} for the vector
$\xi=e_3$), we have the following result.
\begin{state} \label{state:P{3,1a}}
The class $P_{3,1a}$ of potentials that admit the group $G_{3,1a}$
consists of the fields $A_i=A_i(x^4)$.
\end{state}

3.3.1.2. {\it Class $P_{3,1b}$.} %% Point 3.3.1.2
The algebra $\mathcal{L}_{3,1b}=L\{e_1,e_2,e_4\}$ corresponds to the group
$G_{3,1b}$ of translations along the vectors of three-di\-men\-si\-onal
pseudo-Euclidean space $Ox^1x^2x^4$. Since
$\mathcal{L}_{2,1b}\subset\mathcal{L}_{3,1b}$, then the class $C_{3,1b}$
is a subclass of $C_{2,1b}$. Substituting $A_i(x^1,\,x^3)$ for $A_i$ in
the equation \eqref{CondInv A_i:e_4}, we have the following result.
\begin{state} \label{state:P{3,1b}}
The class $P_{3,1b}$ of potentials that admit the group $G_{3,1b}$
consists of the fields $A_i=A_i(x^3)$.
\end{state}

3.3.1.3. {\it Class $P_{3,1c}$} %% Point 3.3.1.3
The algebra $\mathcal{L}_{3,1c}=L\{e_1,e_3,e_2+e_4\}$ corresponds to the group
$G_{3,1c}$ of translations along the vectors of a three-di\-men\-si\-onal
isotropic space. Since $\mathcal{L}_{2,1c}\subset\mathcal{L}_{3,1c}$,
then the class $P_{3,1c}$ is a subclass of $P_{2,1c}$.
Substituting $A_i(x^3,\,x^2-x^4)$ for $A_i$ in \eqref{CondInv A_i:e_3},
we have the following result.
\begin{state} \label{state:P{3,1c}}
The class $P_{3,1c}$ of potentials that admit the group $G_{3,1c}$
consists of the fields $A_i=A_i(x^2-x^4)$.
\end{state}

\subsubsection{Class $P_{3,2}$}  % Subsubsection 3.3.2
The algebra $\mathcal{L}_{3,2}=L\{e_{13}+\lambda e_2,\,e_1,\,e_3\}$
($\lambda\ne0$) corresponds to the group $G_{3,2}$ generated by elliptic
helices with a space-like axis and by translations along the vectors
of the two-di\-men\-si\-onal Euclidian plane. We obtain the class $P_{3,2}$
as a solution of the system of equations \eqref{CondInv A_i:e_1},
\eqref{CondInv A_i:e_3}, and \eqref{CondInv A_i-el} for $\mu=0$:
\begin{equation}\label{CondInv A_i-el(mu=0)}
x^3\partial_1A_i+\lambda\partial_2A_i-x^1\partial_3A_i+
A_1\delta^3_i-A_3\delta^1_i=0.
\end{equation} %\eqno(3.62)
The solution of the system \eqref{CondInv A_i:e_1}--\eqref{CondInv A_i:e_3}
is $A_i=A_i(x^2,\,x^4)$. Substituting $A_i(x^2,\,x^4)$ for $A_i$ in equation
\eqref{CondInv A_i-el(mu=0)}, we have
\begin{equation}\label{P_{3,2}-sist}
\lambda\partial_2A_1-A_3=0,\ \ \lambda\partial_2A_2=0,\ \
\lambda\partial_2A_3+A_1=0,\ \ \lambda\partial_2A_4=0.
\end{equation} %\eqno(3.63)
We obtain the solution of the system \eqref{P_{3,2}-sist} for $\lambda\ne0$
in the form
\begin{equation}\label{Cl:P_{3,2}-ne0}
\begin{split}
&A_1=C_1(x^4)\sin\frac{x^2}{\lambda}+C_2(x^4)\cos\frac{x^2}{\lambda},\ \
 A_2=A_2(x^4),\\
&A_3=C_1(x^4)\cos\frac{x^2}{\lambda}-C_2(x^4)\sin\frac{x^2}{\lambda},\ \
 A_4=A_4(x^4),
\end{split}
\end{equation} %\eqno(3.64)
where $C_1(x^4)$, $C_2(x^4)$, $A_2(x^4)$, and $A_4(x^4)$ are arbitrary
functions.

For $\lambda=0$ the group $G_{3,2}$ is a motion group of the two-dimensional
Euclidian plane; in this case we have the solution of the system
\eqref{P_{3,2}-sist} in the form
\begin{equation}\label{Cl:P_{3,2}-0}
A_1=A_3=0,\ \ A_2=A_2(x^2,\,x^4),\ \ A_4=A_4(x^2,\,x^4).
\end{equation} %\eqno(3.65)
\begin{state} \label{state:P{3,2}}
For $\lambda\ne0$ the class $P_{3,2}$ of potentials that admit
the group $G_{3,2}$ consists of the fields \eqref{Cl:P_{3,2}-ne0};
for $\lambda=0$ this class defined by \eqref{Cl:P_{3,2}-0}.
\end{state}

\subsubsection{Class $P_{3,3}$}  % Subsubsection 3.3.3
The algebra $\mathcal{L}_{3,3}=L\{e_{13}+\mu e_4,\,e_1,\,e_3\}\ (\mu\ne0)$
corresponds to the group $G_{3,3}$ generated by elliptic helices with
a time-like axis and by translations along the vectors
of the two-di\-men\-si\-onal Euclidian plane. We obtain the class $P_{3,3}$
as a solution of the system of equations \eqref{CondInv A_i:e_1},
\eqref{CondInv A_i:e_3}, and \eqref{CondInv A_i-el} for $\lambda=0$:
\begin{equation}\label{CondInv A_i-el(lambda=0)}
x^3\partial_1A_i-x^1\partial_3A_i+\mu\partial_4A_i+
A_1\delta^3_i-A_3\delta^1_i=0.
\end{equation} %\eqno(3.66)
We have the following result.
\begin{state} \label{state:P{3,3}}
The class $P_{3,3}$ of potentials that admit the group $G_{3,3}$
consists of the fields
\begin{equation}\label{Cl:P_{3,3}}
\begin{split}
&A_1=C_1(x^2)\sin\frac{x^4}{\mu}+C_2(x^2)\cos\frac{x^4}{\mu},\ \
 A_2=A_2(x^2),\\
&A_3=C_1(x^2)\cos\frac{x^4}{\mu}-C_2(x^2)\sin\frac{x^4}{\mu},\ \
 A_4=A_4(x^2),
\end{split}
\end{equation} %\eqno(3.67)
where $C_1(x^2)$, $C_2(x^2)$, $A_2(x^2)$, and $A_4(x^2)$ are arbitrary
functions.
\end{state}

\subsubsection{Class $P_{3,4}$}  % Subsubsection 3.3.4
The algebra $\mathcal{L}_{3,4}=L\{e_{13}+\lambda(e_2+e_4),\,e_1,\,e_3\}$
corresponds to the group $G_{3,4}$ generated by elliptic helices
with an isotropic axis and by translations along the vectors
of the two-di\-men\-si\-onal Euclidian plane $Ox^1x^3$. We obtain
the class $P_{3,4}$ as a solution of the system of equations
\eqref{CondInv A_i:e_1}, \eqref{CondInv A_i:e_3} and \eqref{CondInv A_i-el}
for $\lambda=\mu\ne0$:
\begin{equation}\label{CondInv A_i-el(lambda=mu)}
x^3\partial_1A_i+\lambda\partial_2A_i-x^1\partial_3A_i+\lambda\partial_4A_i+
A_1\delta^3_i-A_3\delta^1_i=0.
\end{equation} %\eqno(3.68)
Substituting $A_i(x^2,\,x^4)$ for $A_i$ in \eqref{CondInv A_i-el(lambda=mu)},
we have
\begin{equation}\label{P_{3,4}-sist}
\begin{split}
&\lambda(\partial_2+\partial_4)A_1-A_3=0,\ \
 \lambda(\partial_2+\partial_4)A_2=0,\\
&\lambda(\partial_2+\partial_4)A_3+A_1=0,\ \
 \lambda(\partial_2+\partial_4)A_4=0.
\end{split}
\end{equation} %\eqno(3.69)
Using the substitution
\begin{equation}\label{Change:P_{3,4}}
u=x^2+x^4,\ \ v=x^2-x^4,
\end{equation} %\eqno(3.70)
we get the solution of the system \eqref{P_{3,4}-sist}:
\begin{equation}\label{Cl:P_{3,4}}
\begin{split}
&A_1=C_1(v)\sin\frac{u}{2\lambda}+C_2(v)\cos\frac{u}{2\lambda},\ \
 A_2=A_2(v),\\
&A_3=C_1(v)\cos\frac{u}{2\lambda}-C_2(v)\sin\frac{u}{2\lambda},\ \
 A_4=A_4(v),
\end{split}
\end{equation} %\eqno(3.71)
where $C_1(v)$, $C_2(v)$, $A_2(v)$, and $A_4(v)$ are arbitrary
functions.
\begin{state} \label{state:P{3,4}}
The class $P_{3,4}$ of potentials that admit the group $G_{3,4}$ is defined
by \eqref{Cl:P_{3,4}} and \eqref{Change:P_{3,4}}.
\end{state}

\subsubsection{Class $P_{3,5}$}  % Subsubsection 3.3.5
The algebra $\mathcal{L}_{3,5}=L\{e_{24},\,e_1,\,e_3\}$ corresponds
to the group $G_{3,5}$ generated by pseudo-rotations in the plane $Ox^2x^4$
and by translations along the vectors of the Euclidean plane $Ox^1x^3$.
We obtain the class $P_{3,5}$ as a solution of the system of equations
\eqref{CondInv A_i:e_1}, \eqref{CondInv A_i:e_3}, and \eqref{CondInv A_i-Hyp}
for $\lambda=0$:
\begin{equation}\label{CondInv A_i-Hyp(lambda=0)}
x^4 \partial_2A_i+x^2 \partial_4A_i+A_2\delta^4_i+A_4\delta^2_i=0.
\end{equation} %\eqno(3.72)
Substituting $A_i(x^2,\,x^4)$ for $A_i$ in \eqref{CondInv A_i-Hyp(lambda=0)},
we get
\begin{equation}\label{P_{3,5}-sist}
XA_1=0,\ \ XA_2+A_4=0,\ \ XA_3=0,\ \ XA_4+A_2=0,
\end{equation} %\eqno(3.73)
where $X=x^4\partial_2+x^2\partial_4$. Using the substitution
\begin{equation}\label{Change:P_{3,5}}
x^2=\rho\ch\varphi,\ \ x^4=\rho\sh\varphi,
\end{equation} %\eqno(3.74)
we obtain the solution of the system \eqref{P_{3,5}-sist}in the form
\begin{equation}\label{Cl:P_{3,5}}
\begin{split}
&A_1=A_1(\rho),\ \ A_2=C_1(\rho)\ch\varphi+C_2(\rho)\sh\varphi,\\
&A_3=A_3(\rho),\ \ A_4=-C_1(\rho)\sh\varphi-C_2(\rho)\ch\varphi,
\end{split}
\end{equation} %\eqno(3.75)
where $C_1(\rho)$, $C_2(\rho)$, $A_1(\rho)$, and $A_3(\rho)$ are arbitrary
functions.
\begin{state} \label{state:P{3,5}}
The class $P_{3,5}$ of potentials that admit the group $G_{3,5}$ is defined
by \eqref{Cl:P_{3,5}} and \eqref{Change:P_{3,5}}.
\end{state}

\subsubsection{Class $P_{3,6}$}  % Subsubsection 3.3.6
The algebra $\mathcal{L}_{3,6}=L\{e_{24}+\lambda e_3, e_2, e_4\}$ corresponds
to the group $G_{3,6}$ generated by hyperbolic helices
and translations along the vectors of the pseudo-Euclidean plane.
The algebra $\mathcal{L}_{3,6}$ is an extension of $\mathcal{L}_{2,1b}$
by means of the vector $\xi=e_{24}+\lambda e_3$, therefore
$P_{3,6}\subset P_{2,1b}$. Substituting $A_i(x^1,\,x^3)$ for $A_i$ in
the equation \eqref{CondInv A_i:e_{24}+lambda e_3} (\eqref{CondInv A_i}
for the vector $\xi=e_{24}+\lambda e_3$), we get
\begin{equation}\label{P_{3,6}-sist}
\lambda\partial_3A_1=0,\ \lambda\partial_3A_2+A_4=0,\
\lambda\partial_3A_3=0,\ \lambda\partial_3A_4+A_2=0.
\end{equation} %\eqno(3.76)
For $\lambda\ne0$ we have the following solution of the system
\eqref{P_{3,6}-sist}:
\begin{equation}\label{Cl:P_{3,6}-ne0}
\begin{split}
&A_1=A_1(x^1),\ \
 A_2=C_1(x^1)\ch\frac{x^3}{\lambda}+C_2(x^1)\sh\frac{x^3}{\lambda},\\
&A_3=A_3(x^1),\ \
 A_4=-C_1(x^1)\sh\frac{x^3}{\lambda}-C_2(x^1)\ch\frac{x^3}{\lambda},
\end{split}
\end{equation} %\eqno(3.77)
where $C_1(x^1)$, $C_2(x^1)$, $A_1(x^1)$, and $A_3(x^1)$ are arbitrary
functions.

For $\lambda=0$ $G_{3,6}$ is a motion group of two-dimensional
pseudo--Euclidean plane; in this case we have the solution of the system
\eqref{P_{3,6}-sist} in the form
\begin{equation}\label{Cl:P_{3,6}-0}
A_1=A_1(x^1,\,x^3),\ \ A_2=0,\ \ A_3=A_3(x^1,\,x^3),\ \ A_4=0.
\end{equation} %\eqno(3.78)
\begin{state} \label{state:P{3,6}}
The class $P_{3,6}$ of potentials that admit the group $G_{3,6}$ is defined
by \eqref{Cl:P_{3,6}-ne0} for $\lambda\ne0$; for $\lambda=0$ it is defined by
\eqref{Cl:P_{3,6}-0}.
\end{state}

\subsubsection{Class $P_{3,7}$}  % Subsubsection 3.3.7
The algebra $\mathcal{L}_{3,7}=L\{e_{24}+\lambda e_3,\,e_1,\,e_2-e_4\}$
corresponds to the group $G_{3,7}$ generated by hyperbolic helices
and by translations along the vectors of an isotropic plane.
The algebra $\mathcal{L}_{3,7}$ is an extension of $\mathcal{L}_{2,6}$
by means of the vector $e_1$, therefore $P_{3,7}\subset P_{2,6}$.
Since $A_i$ satisfies to the equation \eqref{CondInv A_i:e_1}, then it is
independent of $x^1=\tilde x^1$; thus we have the following result.
\begin{state}\label{state:P{3,7}}
The class $P_{3,7}$ of potentials that admit the group $G_{3,7}$
consists of the fields
\begin{equation}\label{Cl:P{3,7}}
\begin{split}
&A_1=A_1(u),\ \ A_2=C_1\ch\varphi+C_2\sh\varphi,\\
&A_3=A_3(u),\ \ A_4=-C_1\sh\varphi-C_2\ch\varphi,
\end{split}
\end{equation} %\eqno(3.79)
where
\begin{equation}\label{Cl:P{3,7}-a}
\begin{split}
&C_1=a_1(u)\ch\ln r+a_2(u)\sh\ln r,\\
&C_2=a_1(u)\sh\ln r+a_2(u)\ch\ln r,
\end{split}
\end{equation} %\eqno(3.80)
$u=\tilde{x}^3-\lambda\ln r$, and the transformation of coordinates
is defined by \eqref{Change:Hyperb1,3a-gl2}.
\end{state}

\subsubsection{Class $P_{3,8}$}  % Subsubsection 3.3.8
The algebra
$\mathcal{L}_{3,8}=L\{e_{12}-e_{14}+\lambda e_2$, $e_3$, $e_2-e_4\}$
corresponds to the group $G_{3,8}$ generated by parabolic helices and by
translations along the vectors of an isotropic plane.
The algebra $\mathcal{L}_{3,8}$ is an extension of $\mathcal{L}_{2,7c}$ and
$\mathcal{L}_{2,8}$, therefore $P_{3,8}=P_{2,7c}\cap P_{2,8}$. We have the
following result.
\begin{state} \label{theor:P{3,8}}
The class $P_{3,8}$ of potentials that admit the group $G_{3,8}$
consists of the fields
\begin{equation}\label{Cl:P{3,8}}
\begin{split}
&A_1=C_2\tilde x^2+C_3,\ \ A_2=\frac12C_2(\tilde x^2)^2+
 C_3\tilde x^2+C_1,\\
&A_3=A_3(\tilde x^1),\ \ A_4=A_2+C_2,
\end{split}
\end{equation} %\eqno(3.81)
where $A_3(\tilde x^1)$ and $C_k=C_k(\tilde x^1)$ $(k=1,\,2,\,3)$
are arbitrary functions and $\tilde x^1=2\lambda x^1+\left(x^2+x^4\right)^2$.
\end{state}

\subsubsection{}  % Subsubsection 3.3.9 The class $P_{3,9}$
Here we describe classes of potentials corresponding to the algebra
$\mathcal{L}_{3,9}=L\{e_{12}-e_{14}+\lambda e_2+\mu e_3,\,e_1,\,e_2-e_4\}\ \
(\lambda\mu=0)$ for various $\lambda$ and $\mu$. The corresponding group
$G_{3,9}$ is generated by parabolic helices and by translations along the
vectors of an isotropic plane. The algebra $\mathcal{L}_{3,9}$ is an extension
of $\mathcal{L}_{2,7}$ by means of the vector $e_1$, therefore the
corresponding classes $P_{3,9a}$, $P_{3,9b}$, and $P_{3,9c}$ are restrictions
of classes $P_{2,7a}$, $P_{2,7b}$, and $P_{2,7c}$ by the condition
\eqref{CondInv A_i:e_1}. We consider three cases: a) $\lambda=\mu=0$;
b) $\lambda=0$, $\mu\ne0$; c) $\lambda\ne0$, $\mu=0$.

3.3.9.1. {\it Class $P_{3,9b}$.} %% Point 3.3.9.1
For $\lambda=0$, $\mu\ne0$ we use the substitution \eqref{Change:P{1,4b}},
the equation \eqref{CondInv A_i:e_1} is transformed to
\begin{equation}\label{CondInv A_i:e_1{1,4b}}
-\frac 1{\tilde x^1}\frac{\partial A_i}{\partial\tilde x^2}+
\frac \mu{\tilde x^1}\frac{\partial A_i}{\partial\tilde x^3}-
\tilde x^1\tilde x^2\frac{\partial A_i}{\partial\tilde x^4}=0.
\end{equation} %\eqno(3.82)
Substituting \eqref{Cl:P{2,7a}} for $A_i$ in \eqref{CondInv A_i:e_1{1,4b}},
we get some equation; using a linear independence of functions
$(\tilde x^2)^2$, $\tilde x^2$, and 1, we obtain the following equations
\begin{equation}\label{Sist:P{3,9b}}
\mu\frac{\partial C_2}{\partial\tilde x^3}=0,\ \
\mu\frac{\partial C_3}{\partial\tilde x^3}-C_2=0,\ \
\mu\frac{\partial C_1}{\partial\tilde x^3}-C_3=0,\ \
\mu\frac{\partial A_3}{\partial\tilde x^3}=0
\end{equation} %\eqno(3.83)
for the functions ${C_k(\tilde x^1,\,\tilde x^3)}$ and
${A_3(\tilde x^1,\,\tilde x^3)}$.
We have the solution of \eqref{Sist:P{3,9b}} in the form:
\begin{equation}\label{Cl:P{3,9b}-2}
\begin{split}
&C_1=\frac{(\tilde x^3)^2}{2\mu^2}\Phi(\tilde x^1)+
 \frac{\tilde x^3}{\mu}\Psi(\tilde x^1)+\Xi(\tilde x^1),\ \
 C_2=\Phi(\tilde x^1),\\
&C_3=\frac{\tilde x^3}{\mu}\Phi(\tilde x^1)+\Psi(\tilde x^1),\ \
 A_3=A_3(\tilde x^1),
\end{split}
\end{equation} %\eqno(3.84)
where $A_3(\tilde x^1)$, $\Phi(\tilde x^1)$, $\Psi(\tilde x^1)$ and
$\Xi(\tilde x^1)$ are arbitrary functions and
$$
\tilde x^1=x^2+x^4,\ \ \tilde x^3=x^3+\frac{\mu x^1}{x^2+x^4}.
$$
\begin{state} \label{theor:P{3,9b}}
The class $P_{3,9b}$ of potentials that admit the group $G_{3,9b}$,
corresponding to the algebra $\mathcal{L}_{3,9}$ (${\lambda=0}$,
$\mu\ne0$), consists of the fields
\begin{equation}\label{Cl:P{3,9b}-1}
\begin{split}
&A_1=C_2\tilde x^2+C_3,\ \ A_2=\frac12C_2(\tilde x^2)^2+ C_3\tilde x^2+C_1,\\
&A_3=A_3(\tilde x^1),\ \ A_4=A_2+C_2,
\end{split}
\end{equation} %\eqno(3.85)
where $C_k$ are defined by \eqref{Cl:P{3,9b}-2}.
\end{state}

3.3.9.2. {\it Class $P_{3,9a}$.} %% Point 3.3.9.2
Let now $\mu=\lambda=0$, then $P_{3,9a}\subset P_{2,7a}$; we obtain the
following solution of the system \eqref{Sist:P{3,9b}}
\begin{equation}\label{Cl:P{3,9a}-2}
C_1=C_1(\tilde x^1,\,\tilde x^3),\ \ C_2=C_3=0,\ \
A_3=A_3(\tilde x^1,\,\tilde x^3).
\end{equation} %\eqno(3.86)
\begin{state} \label{theor:P{3,9a}}
The class $P_{3,9a}$ of potentials that admit the group $G_{3,9a}$,
corresponding to the algebra $\mathcal{L}_{3,9}$ ($\lambda=\mu=0$),
consists of the fields
\begin{equation}\label{Cl:P{3,9a}-1}
A_1=0,\ \ A_2=A_4=C_1(\tilde x^1,\,\tilde x^3),\ \
A_3=A_3(\tilde x^1,\,\tilde x^3),
\end{equation} %\eqno(3.87)
where $C_1(\tilde x^1,\,\tilde x^3)$ and $A_3(\tilde x^1,\,\tilde x^3)$
are arbitrary functions and
$$
\tilde x^1=x^2+x^4,\ \ \tilde x^3=x^3.
$$
\end{state}

3.3.9.3. {\it Class $P_{3,9c}$.} %% Point 3.3.9.3
For $\lambda\ne0$, $\mu=0$ the algebra
$$
\mathcal{L}_{3,9c}=L\{e_{12}-e_{14}+\lambda e_2,\,e_1,\,e_2-e_4\}
$$
includes the algebra $\mathcal{L}_{2,7c}$, therefore
$P_{3,9c}\subset P_{2,7c}$. Here we use the substitution
\eqref{Change:P{1,4c}}, the equation \eqref{CondInv A_i:e_1}
takes the form
\begin{equation}\label{CondInv A_i:e_1{1,4c}}
2\lambda\frac{\partial A_i}{\partial\tilde x^1}+
\lambda\tilde x^2\frac{\partial A_i}{\partial\tilde x^4}=0.
\end{equation} %\eqno(3.88)
Since $A_i$ is independent of $\tilde x^4$, we have the following result.
\begin{state} \label{theor:P{3,9c}}
The class $P_{3,9c}$ of potentials that admit the group $G_{3,9c}$,
corresponding to the algebra $\mathcal{L}_{3,9c}$, consists of the fields
\begin{equation}\label{Cl:P{3,9c}}
\begin{split}
&A_1=C_2\tilde x^2+C_3,\ \ A_2=\frac12C_2(\tilde x^2)^2+ C_3\tilde x^2+C_1,\\
&A_3=A_3(\tilde x^3),\phantom{===} A_4=A_2+C_2,
\end{split}
\end{equation} %\eqno(3.89)
where $A_3(\tilde x^3)$ and $C_k=C_k(\tilde x^3)$ $(k=1,\,2,\,3)$
are arbitrary functions and
$\tilde x^2=(x^2+x^4)/{\lambda},\ \ \tilde x^3=x^3$.
\end{state}

\subsubsection{}  % Subsubsection 3.3.10 The class $P_{3,10}$
Here we describe classes of potentials corresponding to the algebra
$\mathcal{L}_{3,10}=L\{e_{12}-e_{14}+\lambda e_2,\,e_1+\mu e_3,\,e_2-e_4\}$
for various $\lambda$ and $\mu$. The corresponding group $G_{3,10}$
is generated by parabolic helices or parabolic rotations and by translations
along the vectors of an isotropic plane. If $\mu=0$, then
$\mathcal{L}_{3,10}=\mathcal{L}_{3,9c}$. We consider two cases:
a) $\lambda\ne0$, $\mu\ne0$; b) $\lambda=0$, $\mu\ne0$.

3.3.10.1. {\it Class $P_{3,10a}$.} % Point 3.10.1.1
Let $\lambda\ne0$, $\mu\ne0$. The algebra
$\mathcal{L}_{3,10a}=\mathcal{L}_{3,10}$ is an extension
of $\mathcal{L}_{2,7c}$ by means of the vector $e_1+\mu e_3$, therefore
$P_{3,10a}\subset P_{2,7c}$. The equation \eqref{CondInv A_i}
for ${\xi=e_1+\mu e_3}$ takes the form
\begin{equation}\label{CondInv A_i:e_1+mu e_3}
\partial_1A_i+\mu\partial_3A_i=0.
\end{equation} %\eqno(3.90)
Using the substitution \eqref{Change:P{1,4c}} we transform
\eqref{CondInv A_i:e_1+mu e_3} to the form
\begin{equation}\label{CondInv A_i:e_1+mu e_3(1,4c)}
2\lambda\frac{\partial A_i}{\partial\tilde x^1}+
\mu\frac{\partial A_i}{\partial\tilde x^3}+
\lambda\tilde x^2\frac{\partial A_i}{\partial\tilde x^4}=0.
\end{equation} %\eqno(3.91)
Substituting \eqref{Cl:P{2,7a}}--\eqref{Change:P{1,4c}} for $A_i$ in
\eqref{CondInv A_i:e_1+mu e_3(1,4c)}, we get the following result.
\begin{state} \label{theor:P{3,10a}}
The class $P_{3,10a}$ of potentials that admit the group $G_{3,10a}$,
corresponding to the algebra $\mathcal{L}_{3,10a}$, consists of the fields
\begin{equation}\label{Cl:P{3,10a}}
\begin{split}
&A_1=C_2\tilde x^2+C_3,\ \ A_2=\frac12C_2(\tilde x^2)^2+ C_3\tilde x^2+C_1,\\
&A_3=A_3(\mu\tilde x^1-2\lambda\tilde x^3),\ \ A_4=A_2+C_2,
\end{split}
\end{equation} %\eqno(3.92)
where $A_3(\mu\tilde x^1-2\lambda\tilde x^3)$ and
$C_k=C_k(\mu\tilde x^1-2\lambda\tilde x^3)$ $(k=1,\,2,\,3)$
are arbitrary functions and
$$
\tilde x^1=2\lambda x^1+\left(x^2+x^4\right)^2,\ \ \tilde x^3=x^3.
$$
\end{state}

3.3.10.2. {\it Class $P_{3,10b}$} %% Point 3.3.10.2
Let $\lambda=0$, $\mu\ne0$. The algebra
$$
\mathcal{L}_{3,10b}=L\{e_{12}-e_{14},\,e_1+\mu e_3,\,e_2-e_4\}
$$
is an extension of $\mathcal{L}_{2,7a}$ by means of the vector $e_1+\mu e_3$,
hence $P_{3,10b}\subset P_{2,7a}$. By means of substitution
\eqref{Change:P{1,4a}} equation \eqref{CondInv A_i:e_1+mu e_3} is
transformed to the form
\begin{equation}\label{CondInv A_i:e_1+mu e_3(1,4a)}
-\frac1{\tilde x^1}\frac{\partial A_i}{\partial\tilde x^2}+
\mu\frac{\partial A_i}{\partial\tilde x^3}=0.
\end{equation} %\eqno(3.93)
Substituting \eqref{Cl:P{2,7a}} for $A_i$ in
\eqref{CondInv A_i:e_1+mu e_3(1,4a)}, we get some equation;
using a linear independence of functions $(\tilde x^2)^2$, $\tilde x^2$,
and 1, we obtain the following equations
\begin{equation}\label{Sist:P{3,10b}}
\mu\frac{\partial C_2}{\partial\tilde x^3}=0,\ \
\mu\frac{\partial C_3}{\partial\tilde x^3}-\frac{C_2}{\tilde x^1}=0,\ \
\mu\frac{\partial C_1}{\partial\tilde x^3}-\frac{C_3}{\tilde x^1}=0,\ \
\mu\frac{\partial A_3}{\partial\tilde x^3}=0
\end{equation} %\eqno(3.94)
for the functions ${C_k(\tilde x^1,\,\tilde x^3)}$ and
${A_3(\tilde x^1,\,\tilde x^3)}$.
We have the solution of \eqref{Sist:P{3,10b}} in the form:
\begin{equation}\label{Cl:P{3,10b}}
\begin{split}
&C_1=\frac{(\tilde x^3)^2}{2\mu^2(\tilde x^1)^2}\Phi(\tilde x^1)+
 \frac{\tilde x^3}{\mu\tilde x^1}\Psi(\tilde x^1)+\Xi(\tilde x^1),\ \
 C_2=\Phi(\tilde x^1),\\
&C_3=\frac{\tilde x^3}{\mu\tilde x^1}\Phi(\tilde x^1)+\Psi(\tilde x^1),\ \
 A_3=A_3(\tilde x^1),
\end{split}
\end{equation} %\eqno(3.95)
where $A_3(\tilde x^1)$, $\Phi(\tilde x^1)$, $\Psi(\tilde x^1)$, and
$\Xi(\tilde x^1)$ are arbitrary functions.
\begin{state} \label{theor:P{3,10b}}
The class $P_{3,10b}$ of potentials that admit the group $G_{3,10b}$,
corresponding to the algebra $\mathcal{L}_{3,10b}$, consists of the fields
\eqref{Cl:P{3,9b}-1}, where $C_k$ are defined by \eqref{Cl:P{3,10b}} and
$$
\tilde x^1=x^2+x^4,\ \ \tilde x^2=-\frac{x^1}{x^2+x^4},\ \ \tilde x^3=x^3.
$$
\end{state}

\subsubsection{Class $P_{3,11}$}  % Subsubsection 3.3.11
The algebra $\mathcal{L}_{3,11}=L\{e_{13}+\lambda e_{24},\,e_1,\,e_3\}$
corresponds to the group $G_{3,11}$ generated by proportional bi-rotations
and by translations along the vectors of two-di\-men\-si\-onal Euclidian plane.
The class $P_{3,11}$, corresponding to the algebra $\mathcal{L}_{3,11}$,
is a subclass of the class $P_{1,5}$. For description of it we substitute
\eqref{Cl:P{1,5}} for $A_i$ in equations
\begin{equation}\label{CondInv A_i:e_1-ch}
\cos(\theta-\varphi)\frac{\partial A_i}{\partial r}-
\frac{\sin(\theta-\varphi)}r\frac{\partial A_i}{\partial\theta}=0
\end{equation} %% \eqno(3.96)
and
\begin{equation}\label{CondInv A_i:e_3-ch}
\sin(\theta-\varphi)\frac{\partial A_i}{\partial r}+
\frac{\cos(\theta-\varphi)}r\frac{\partial A_i}{\partial\theta}=0
\end{equation} %% \eqno(3.97)
(equations \eqref{CondInv A_i:e_1} and \eqref{CondInv A_i:e_3}, transformed
by substitution \eqref{Change:P{1,5}}); we have the following solution
\begin{equation}\label{Cl:P{3,11}}
\begin{split}
&A_1=C_1(\rho)\cos\varphi+C_2(\rho)\sin\varphi,\\
&A_2=C_3(\rho)\ch\lambda\varphi+C_4(\rho)\sh\lambda\varphi,\\
&A_3=-C_1(\rho)\sin\varphi+C_2(\rho)\cos\varphi,\\
&A_4=-C_3(\rho)\sh\lambda\varphi-C_4(\rho)\ch\lambda\varphi,
\end{split}
\end{equation} %% \eqno(3.98)
where $C_k=C_k(\rho)$ are arbitrary functions.
\begin{state} \label{theor:P{3,11}}
The class $P_{3,11}$ of potentials that admit the group $G_{3,11}$
consists of the fields \eqref{Cl:P{3,11}}.
\end{state}

\subsubsection{Class $P_{3,12}$}  % Subsubsection 3.3.12
The algebra $\mathcal{L}_{3,12}=L\{e_{13}+\lambda e_{24},\,e_2,\,e_4\}$
corresponds to the group $G_{3,12}$ generated by proportional bi-rotations
and by translations along the vectors of two-di\-men\-si\-onal
pseudo-Euclidian plane. The algebra $\mathcal{L}_{3,12}$ is an extension
of $\mathcal{L}_{1,5}$ by means of the vectors $e_2$ and $e_4$, therefore
$P_{3,12}\subset P_{1,5}$. For description of it we substitute
\eqref{Cl:P{1,5}} for $A_i$ in equations
\begin{equation}\label{CondInv A_i:e_2-ch}
\ch(\lambda\varphi)\frac{\partial A_i}{\partial\rho}-
\frac{\sh(\lambda\varphi)}{\lambda\rho}\frac{\partial A_i}{\partial\theta}-
\frac{\sh(\lambda\varphi)}{\lambda\rho}\frac{\partial A_i}{\partial\varphi}=0
\end{equation} %% \eqno(3.99)
and
\begin{equation}\label{CondInv A_i:e_4-ch}
-\sh(\lambda\varphi)\frac{\partial A_i}{\partial\rho}+
\frac{\ch(\lambda\varphi)}{\lambda\rho}\frac{\partial A_i}{\partial\theta}+
\frac{\ch(\lambda\varphi)}{\lambda\rho}\frac{\partial A_i}{\partial\varphi}=0
\end{equation} %% \eqno(3.100)
(equations \eqref{CondInv A_i:e_2} and \eqref{CondInv A_i:e_4}, transformed
by substitution \eqref{Change:P{1,5}}). Multiplying \eqref{CondInv A_i:e_2-ch}
by $\ch(\lambda\varphi)$, \eqref{CondInv A_i:e_4-ch} by $\sh(\lambda\varphi)$
and summing received equations, we get
\begin{equation}\label{Eqv:P{3,12}-1}
\partial A_i/\partial\rho=0;
\end{equation} %% \eqno(3.101)
therefore $A_i$ and $C_k$ are independent of $\rho$: $C_k=C_k(r,\theta)$.
Further, we have the following consequence of equations
\eqref{CondInv A_i:e_2-ch} and \eqref{Eqv:P{3,12}-1}:
\begin{equation}\label{Eqv:P{3,12}-2}
\frac{\partial A_i}{\partial\theta}+\frac{\partial A_i}{\partial\varphi}=0.
\end{equation} %% \eqno(3.102)
The system \eqref{CondInv A_i:e_2-ch}--\eqref{CondInv A_i:e_4-ch} is
equivalent to the system \eqref{Eqv:P{3,12}-1}--\eqref{Eqv:P{3,12}-2}.
Substituting \eqref{Cl:P{1,5}} for $A_i$ in
\eqref{Eqv:P{3,12}-1}--\eqref{Eqv:P{3,12}-2}, we get a result
of calculations:
\begin{equation}\label{Cl:P{3,12}}
\begin{split}
&A_1=a_1(r)\sin(\theta-\varphi)+a_2(r)\cos(\theta-\varphi),\\
&A_2=a_3(r)\sh[\lambda(\theta-\varphi)]+a_4(r)\ch[\lambda(\theta-\varphi)],\\
&A_3=-a_1(r)\cos(\theta-\varphi)+a_2(r)\sin(\theta-\varphi),\\
&A_4=a_3(r)\ch[\lambda(\theta-\varphi)]+a_4(r)\sh[\lambda(\theta-\varphi)],
\end{split}
\end{equation} %% \eqno(3.103)
where $a_k=a_k(r)$ are arbitrary functions.
\begin{state} \label{theor:P{3,12}}
The class $P_{3,12}$ of potentials that admit the group $G_{3,12}$
consists of the fields \eqref{Cl:P{3,12}} (the transformation of coordinates
is defined by \eqref{Change:P{1,5}}).
\end{state}

\subsubsection{Class $P_{3,13}$}  % Subsubsection 3.3.13
The algebra $\mathcal{L}_{3,13}=L\{e_{13},\,e_{24},\,e_2-e_4\}$
is an extension of the algebra $\mathcal{L}_{2,10}$ by means of the vector
$e_2-e_4$, therefore $P_{3,13}\subset P_{2,10}$. By substitution
\eqref{Change:P{2,10}} the equation \eqref{CondInv A_i:e_2-e_4}
(\eqref{CondInv A_i} for the vector $\xi=e_2-e_4$) takes the form
\begin{equation}\label{Eqv:P{3,13}}
\frac{\partial A_i}{\partial\varphi}+\frac{\partial A_i}{\partial\theta}-
\rho\frac{\partial A_i}{\partial\rho}=0.
\end{equation} %% \eqno(3.104)
Substituting \eqref{Cl:P{2,10}} for $A_i$ in the equation \eqref{Eqv:P{3,13}},
we obtain the following result.
\begin{state} \label{theor:P{3,13}}
The class $P_{3,13}$ of potentials $A_i$ that admit the group $G_{3,13}$
consists of the fields
\begin{equation}\label{Cl:P{2,13}}
\begin{split}
&A_1=-t_1(r)\sin(\theta-\varphi)+t_2(r)\cos(\theta-\varphi),\\
&A_2=\rho\,C(r)e^{\varphi}+\frac{D(r)}{\rho}e^{-\varphi},\\
&A_3=t_1(r)\cos(\theta-\varphi)+t_2(r)\sin(\theta-\varphi),\\
&A_4=-\rho\,C(r)e^{\varphi}+\frac{D(r)}{\rho}e^{-\varphi},
\end{split}
\end{equation} %\eqno(3.105)
where $t_1(r)$, $t_2(r)$, $C(r)$, and $D(r)$ are arbitrary functions.
(the transformation of coordinates is defined by \eqref{Change:P{2,10}}).
\end{state}

\subsubsection{Class $P_{3,14}$}  % Subsubsection 3.3.14
The algebra
$$
\mathcal L_{3,14}=L\{e_{12}-e_{14}+\lambda e_1+\mu e_3,\,
e_{23}+e_{34}+\nu e_1+\lambda e_3,\,e_2-e_4\}
$$
corresponds to the group $G_{3,14}$ generated by two one-di\-men\-si\-onal
subgroups of parabolic helices and by translations along an
isotropic straight line. The equation \eqref{CondInv A_i} for basis vectors
of the algebra $\mathcal L_{3,14}$ take the following forms
\begin{gather}
XA_i-A_1(\delta^2_i+\delta^4_i)+(A_2-A_4)\delta^1_i=0,\nonumber\\
 X=(\lambda-x^2-x^4)\partial_1+\mu\partial_3,
\label{CondInv A_i:e_{12}-e_{14}+lambda e_1+mu e_3}\\ %\eqno(3.106)
YA_i-(A_2-A_4)\delta^3_i+A_3(\delta^2_i+\delta^4_i)=0,\nonumber\\
 Y=\nu\partial_1+(\lambda+x^2+x^4)\partial_3,
\label{CondInv A_i:e_{23}+e_{34}+nu e_1+lambda e_3} %\eqno(3.107)
\end{gather}
and \eqref{CondInv A_i:e_2-e_4}. We have the solution of the equation
\eqref{CondInv A_i:e_2-e_4} in the form
\begin{equation}\label{A_i(e_2-e_4)}
A_i=A_i(x^1,\,x^2+x^4,\,x^3).
\end{equation} %\eqno(3.108)
Substituting \eqref{A_i(e_2-e_4)} for $A_i$ in equations
\eqref{CondInv A_i:e_{12}-e_{14}+lambda e_1+mu e_3} and
\eqref{CondInv A_i:e_{23}+e_{34}+nu e_1+lambda e_3} and transforming theirs
by means the substitution
\begin{equation}\label{Change:P{3,14}}
u=x^2+x^4,\ \varphi=\frac{\mu x^1+(u-\lambda)x^3}{u^2-\lambda^2+\mu\nu},\
\psi=\frac{\nu x^3-(u+\lambda)x^1}{u^2-\lambda^2+\mu\nu},
\end{equation} %\eqno(3.109)
we obtain two systems
\begin{align}\label{CondInv A_i:e_{12}-e_{14}+lambda e_1+mu e_3-ch}
&\frac{\partial A_1}{\partial\psi}+A_2-A_4=0,\,
 \frac{\partial A_2}{\partial\psi}-A_1=0,\,
 \frac{\partial A_3}{\partial\psi}=0,\,
 \frac{\partial A_4}{\partial\psi}-A_1=0,\\ %\eqno(3.110)
              \label{CondInv A_i:e_{23}+e_{34}+nu e_1+lambda e_3-ch}
&\frac{\partial A_1}{\partial\varphi}=0,\,
 \frac{\partial A_2}{\partial\varphi}+A_3=0,\,
 \frac{\partial A_3}{\partial\varphi}-A_2+A_4=0,\,
 \frac{\partial A_4}{\partial\varphi}+A_3=0. %\eqno(3.111)
\end{align}
We have the total solution of the system
\eqref{CondInv A_i:e_{12}-e_{14}+lambda e_1+mu e_3-ch} in the form
\begin{equation}\label{Sol A_i:e_{12}-e_{14}+lambda e_1+mu e_3-ch}
\begin{split}
&A_1=\psi\,C_2(u,\varphi)+C_3(u,\varphi),\\
&A_2=\frac12\psi^2\,C_2(u,\varphi)+\psi\,C_3(u,\varphi)+C_1(u,\varphi),\\
&A_3=A_3(u,\varphi),\ \ A_4=A_2+C_2(u,\varphi).
\end{split}
\end{equation} %\eqno(3.112)
Finally, substituting \eqref{Sol A_i:e_{12}-e_{14}+lambda e_1+mu e_3-ch}
for $A_i$ in \eqref{CondInv A_i:e_{23}+e_{34}+nu e_1+lambda e_3-ch}, we obtain
\begin{equation}\label{Cl:P{3,14}}
A_1=C_3(u),\ \ A_2=A_4=\psi\,C_3(u)+C_1(u),\ \ A_3=0,
\end{equation} %\eqno(3.113)
where $C_1(u)$ and $C_3(u)$ are arbitrary functions.
\begin{state} \label{theor:P{3,14}}
The class $P_{3,14}$ of potentials that admit the group $G_{3,14}$
consists of the fields defined by \eqref{Cl:P{3,14}} and
\eqref{Change:P{3,14}}.
\end{state}

\subsubsection{Class $P_{3,15}$}  % Subsubsection 3.3.15
The algebra $\mathcal{L}_{3,15}=L\{e_{12}-e_{14},\ e_{24},\ e_3\}$
corresponds to the group $G_{3,15}$ generated by parabolic rotations, by
pseudo-rotations, and by translations along a space-like straight line.
The algebra $\mathcal{L}_{3,15}$ is an extension of the algebra
$\mathcal{L}_{2,12a}$ ($\mathcal{L}_{2,12}$ for $\lambda=0$) by means of the
vector $e_3$, therefore the class $P_{3,15}$ is a subclass of $P_{2,12a}$
($P_{2,12}$ for $\lambda=0$). Substituting
\eqref{Cl:P_{2,12}}--\eqref{Change:P_{2,12}}--\eqref{Change:P{1,4a}} for $A_i$
in equation \eqref{CondInv A_i:e_3}, which means independence $A_i$ of
$\tilde x^3$, we obtain the following result.
\begin{state} \label{theor:P{3,15}}
The class $P_{3,15}$ of potentials that admit the group $G_{3,15}$
consists of the fields defined by \eqref{Cl:P_{2,12}}, where
$$
\Phi_k=
\Phi_k(v)=\Phi_k\left(\tilde x^4-\frac12\left(\tilde x^1\right)^2\right)
$$
are arbitrary functions and the transformation of coordinates is defined
by \eqref{Change:P{1,4a}}.
\end{state}

\subsubsection{Class $P_{3,16}$}  %% Subsubsection 3.3.16
The algebra
$$
\mathcal{L}_{3,16}=L\{e_{12}-e_{14},\,e_{24}+\lambda e_1+\mu e_3,\,e_2-e_4\}
$$
corresponds to the group $G_{3,16}$ generated by parabolic rotations, by
hyperbolic helices, and by translations along an isotropic straight line.
The algebra $\mathcal{L}_{3,16}$ is an extension of the algebra
$\mathcal{L}_{2,7a}=L\{{e_{12}-e_{14}},\break {e_2-e_4}\}$ by means of the
vector $e_{24}+\lambda e_1+\mu e_3$, therefore $P_{3,16}\subset P_{2,7a}$.
The equation \eqref{CondInv A_i} for the vector
$\xi=e_{24}+\lambda e_1+\mu e_3$ takes the form
\begin{gather}\label{CondInv A_i:e_{24}+lambda e_1+mu e_3}
XA_i+A_2\delta^4_i+A_4\delta^2_i=0,\\ %\eqno(3.114)
              \label{X:e_{24}+lambda e_1+mu e_3-ch}
Xf=\lambda\partial_1f+x^4\partial_2f+\mu\partial_3f+x^2\partial_4f=
  -\lambda\left(\frac 1{\tilde x^1}\frac{\partial f}{\partial\tilde x^2}+
  \tilde x^1\tilde x^2\frac{\partial f}{\partial\tilde x^4}\right)+\nonumber\\
  +\tilde x^1\frac{\partial f}{\partial\tilde x^1}-
\tilde x^2\frac{\partial f}{\partial\tilde x^2}+
\mu\frac{\partial f}{\partial\tilde x^3}+
(\tilde x^1)^2\frac{\partial f}{\partial\tilde x^4} %\eqno(3.115)
\end{gather}
(here we use the substitution \eqref{Change:P{1,4a}}). Substituting
\eqref{Cl:P{2,7a}} for $A_i$ in \eqref{CondInv A_i:e_{24}+lambda e_1+mu e_3},
we get some equation; using a linear independence of functions
$(\tilde x^2)^2$, $\tilde x^2$, and 1, we obtain the following equations
\begin{equation}\label{Equation:P_{3,16}}
\begin{split}
&\tilde x^1\frac{\partial C_2}{\partial\tilde x^1}+
 \mu\frac{\partial C_2}{\partial\tilde x^3}-C_2=0,\ \
 \tilde x^1\frac{\partial C_3}{\partial\tilde x^1}+
 \mu\frac{\partial C_3}{\partial\tilde x^3}-\frac{\lambda}{\tilde x^1}C_2=0,\\
&\tilde x^1\frac{\partial C_1}{\partial\tilde x^1}+
 \mu\frac{\partial C_1}{\partial\tilde x^3}+C_1+C_2-
 \frac{\lambda}{\tilde x^1}C_3=0,\ \
 \tilde x^1\frac{\partial A_3}{\partial\tilde x^1}+
 \mu\frac{\partial A_3}{\partial\tilde x^3}=0
\end{split}
\end{equation} %\eqno(3.116)
for the functions $C_k(\tilde x^1,\tilde x^3)$ and
$A_3(\tilde x^1,\tilde x^3)$. Using the substitution
\begin{equation}\label{Change:P_{3,16}}
u=\tilde x^3-\mu\ln\tilde x^1,\ \ v=\ln\tilde x^1,
\end{equation} %\eqno(3.117)
we transform \eqref{Equation:P_{3,16}}; the solution of received system
takes the form
\begin{equation}\label{Cl:P_{3,16}-2}
\begin{split}
&C_1=\Phi_3(u)\,e^{-v}-\frac12\Phi_1(u)\,e^{v}+
 v\,e^{-v}\left[\frac12\lambda^2v\Phi_1(u)+\lambda\Phi_2(u)\right],\\
&C_2=\Phi_1(u)\,e^{v},\ \ C_3=\lambda\,v\,\Phi_1(u)+\Phi_2(u),\ \
 A_3=\Phi_4(u),
\end{split}
\end{equation} %\eqno(3.118)
where $\Phi_k(u)$ are arbitrary functions. Thus we have the following result.
\begin{state} \label{theor:P{3,16}}
The class $P_{3,16}$ of potentials that admit the group $G_{3,16}$
consists of the fields
\begin{equation}\label{Cl:P{3,16}-1}
\begin{split}
&A_1=C_2\tilde x^2+C_3,\ \ A_2=\frac12C_2(\tilde x^2)^2+ C_3\tilde x^2+C_1,\\
&A_3=\Phi_4(u),\ \ A_4=A_2+C_2,
\end{split}
\end{equation} %\eqno(3.119)
where $C_l=C_l(u,v)$ are defined by
\eqref{Cl:P_{3,16}-2}--\eqref{Change:P_{3,16}}, $\Phi_k(u)$
are arbitrary functions, and
$$
\tilde x^1=x^2+x^4,\ \ \tilde x^2=-\frac{x^1}{x^2+x^4},\ \ \tilde x^3=x^3.
$$
\end{state}

\subsubsection{Class $P_{3,17}$}  % subsubsection 3.3.17
The algebra $\mathcal{L}_{3,17}=L\{e_{12}-e_{14},\,e_{23}+e_{34},\,e_{24}\}$
corresponds to the group $G_{3,17}$ generated by two one-di\-men\-si\-onal
subgroups of parabolic helices and by hyperbolic rotations.
The algebra $\mathcal{L}_{3,17}$ is an extension of the algebra
$\mathcal{L}_{2,11a}$ by means of the vector $e_{24}$, hence
$P_{3,17}\subset P_{2,11a}$. The equation \eqref{CondInv A_i} for the vector
$\xi=e_{24}$ takes the form
\begin{equation}\label{CondInv A_i:e_{24}}
XA_1=0,\ \ XA_2+A_4=0,\ \ XA_3=0,\ \ XA_4+A_2=0,
\end{equation} %\eqno(3.120)
where $X=x^4\partial_2+x^2\partial_4$.
By substitution \eqref{Change:P_{2,11a}} we replace $X$ to the form
\begin{equation}\label{X:e_{24}-ch}
X=\tilde x^1\frac{\partial}{\partial\tilde x^1}-
 \tilde x^2\frac{\partial}{\partial\tilde x^2}-
 \tilde x^3\frac{\partial}{\partial\tilde x^3}.
\end{equation} %\eqno(3.121)
Substituting \eqref{Solution:P_{2,11a}} for $A_i$ in
\eqref{CondInv A_i:e_{24}}--\eqref{X:e_{24}-ch}, we get some equation;
using a linear independence of functions $(\tilde x^2)^2$, $\tilde x^2$,
$(\tilde x^3)^2$, $\tilde x^3$ and 1, we obtain some system of equations;
solving this system, we get the following result.
\begin{state} \label{state:P_{3,17}}
The class $P_{3,17}$ of potentials that admit the group $G_{3,17}$
consists of the fields defined by \eqref{Solution:P_{2,11a}}, where
\begin{equation}\label{Phi,Psi,Xi,Theta:P_{3,17}}
\begin{split}
&\Phi=\tilde x^1C_1(\tilde x^4),\ \ \Psi=C_2(\tilde x^4),\ \
\Xi=C_3(\tilde x^4),\\
&\Theta=\frac{\tilde x^1}2C_1(\tilde x^4)+ \frac{C_4(\tilde x^4)}{\tilde x^1},
\end{split}
\end{equation} %\eqno(3.122)
$C_k(\tilde x^4)$ are arbitrary functions, and the transformation
of coordinates is defined by \eqref{Change:P_{2,11a}}.
\end{state}

\subsubsection{}  % subsubsection 3.3.18 Class $P_{3,18}$
The algebra $\mathcal{L}_{3,18}=L\{e_{12}-e_{14},\,e_{23}+e_{34},\,
e_{13}+\lambda(e_2-e_4)\}$
corresponds to the group $G_{3,18}$ generated by two one-di\-men\-si\-onal
subgroups of parabolic rotations and by elliptic helices with an isotropic
axis (or rotations for ${\lambda=0}$). The algebra $\mathcal{L}_{3,18}$
is an extension of the algebra $\mathcal{L}_{2,11a}$ by means of the vector
$e_{13}+\lambda(e_2-e_4)$, hence $P_{3,18}\subset P_{2,11a}$.
The equation \eqref{CondInv A_i} for the vector $\xi=e_{13}+\lambda(e_2-e_4)$
takes the form
\begin{equation}\label{CondInv A_i:e_{13}+lambda(e_2-e_4)-3,18}
XA_1-A_3=0,\ \ XA_2=0,\ \ XA_3+A_1=0,\ \ XA_4=0,
\end{equation} %\eqno(3.123)
where $X=x^3\partial_1+\lambda\partial_2-x^1\partial_3-\lambda\partial_4$.
By substitution \eqref{Change:P_{2,11a}} we replace $X$ to the form
\begin{equation}\label{X:e_{13}+lambda(e_2-e_4)-ch}
X=-\tilde x^3\frac{\partial}{\partial\tilde x^2}+
 \tilde x^2\frac{\partial}{\partial\tilde x^3}+
 2\lambda\tilde x^1\frac{\partial}{\partial\tilde x^4}.
\end{equation} %\eqno(3.124)
Substituting \eqref{Solution:P_{2,11a}} for $A_i$ in
\eqref{CondInv A_i:e_{13}+lambda(e_2-e_4)-3,18}--%
\eqref{X:e_{13}+lambda(e_2-e_4)-ch}, we get some equation;
using a linear independence of variables $\tilde x^2$, $\tilde x^3$, and
their powers, we obtain the system of equations:
\begin{equation}\label{Equations:P{3,18}}
\lambda\frac{\partial\Phi}{\partial\tilde x^4}=0,\
2\lambda\tilde x^1\frac{\partial\Psi}{\partial\tilde x^4}-\Xi=0,\
\lambda\frac{\partial\Theta}{\partial\tilde x^4}=0,\
2\lambda\tilde x^1\frac{\partial\Xi}{\partial\tilde x^4}+\Psi=0.
\end{equation} %\eqno(3.125)

3.3.18.1. {\it Class $P_{3,18a}$.} %% Point 3.3.18.1
For ${\lambda\ne0}$ the solution of equations \eqref{Equations:P{3,18}}
takes the form
\begin{equation}\label{Phi,Psi,Xi,Theta:P{3,18a}}
\begin{split}
&\Psi=C_1(\tilde x^1)\cos\frac{\tilde x^4}{2\lambda\tilde x^1}+
 C_2(\tilde x^1)\sin\frac{\tilde x^4}{2\lambda\tilde x^1},\ \
 \Phi=C_3(\tilde x^1),\\
&\Xi=-C_1(\tilde x^1)\sin\frac{\tilde x^4}{2\lambda\tilde x^1}+
 C_2(\tilde x^1)\cos\frac{\tilde x^4}{2\lambda\tilde x^1},\ \
 \Theta=C_4(\tilde x^1),
\end{split}
\end{equation} %\eqno(3.126)
where $C_k(\tilde x^1)$ are arbitrary functions.
\begin{state} \label{state:P_{3,18a}}
The class $P_{3,18a}$ of potentials that admit the group $G_{3,18a}$
corresponding to the algebra $\mathcal{L}_{3,18}$ (${\lambda\ne0}$)
is defined by \eqref{Solution:P_{2,11a}} and \eqref{Phi,Psi,Xi,Theta:P{3,18a}}
(the transformation of coordinates is defined by \eqref{Change:P_{2,11a}}).
\end{state}

3.3.18.2. {\it Class $P_{3,18b}$.} %% Point 3.3.18.2
For ${\lambda=0}$ the solution of equations \eqref{Equations:P{3,18}}
takes the form
\begin{equation}\label{Phi,Psi,Xi,Theta:P{3,18b}}
\Phi=\Phi(\tilde x^1,\tilde x^4),\ \ \Psi=\Xi=0,\ \
\Theta=\Theta(\tilde x^1,\tilde x^4),
\end{equation} %\eqno(3.127)
where ${\Phi(\tilde x^1,\tilde x^4)}$ and $\Theta(\tilde x^1,\tilde x^4)$
are arbitrary functions. Substituting \eqref{Phi,Psi,Xi,Theta:P{3,18b}}
for $\Phi$, $\Psi$, $\Xi$, and $\Theta$ in \eqref{Solution:P_{2,11a}},
we get the following result.
\begin{state} \label{state:P_{3,18b}}
The class $P_{3,18b}$ of potentials that admit the group $G_{3,18b}$
corresponding to the algebra
$\mathcal{L}_{3,18b}=L\{e_{12}-e_{14},\,e_{23}+e_{34},\,e_{13}\}$
($\mathcal{L}_{3,18}$ for ${\lambda=0}$) is defined by the following formulae
\begin{align}\label{Cl:P_{3,18b}}
&A_1=-\tilde x^2\,\Phi(\tilde x^1,\tilde x^4),\ \
 A_2=-\frac12((\tilde x^2)^2+(\tilde x^3)^2)\Phi(\tilde x^1,\tilde x^4)+
 \Theta(\tilde x^1,\tilde x^4),\nonumber\\
&A_3=\tilde x^3\,\Phi(\tilde x^1,\tilde x^4),\ \
 A_4=A_2-\Phi(\tilde x^1,\tilde x^4)
\end{align} %\eqno(3.128)
(the transformation of coordinates is defined by \eqref{Change:P_{2,11a}}).
\end{state}

\subsubsection{Class $P_{3,19}$}  % subsubsection 3.3.19
The algebra
$$
\mathcal{L}_{3,19}=L\{e_{12}-e_{14},\,e_{23}+e_{34},\,e_{13}+\lambda e_{24}\}\
\ (\lambda\ne0)
$$
corresponds to the group $G_{3,19}$ generated by two one-di\-men\-si\-onal
subgroups of parabolic rotations and by one-di\-men\-si\-onal subgroup of
be-rotations. The algebra $\mathcal{L}_{3,19}$ is an extension of the algebra
$\mathcal{L}_{2,11a}$ by means of the vector $e_{13}+\lambda e_{24}$, hence
$P_{3,19}\subset P_{2,11a}$. The equation \eqref{CondInv A_i} for the vector
$\xi=e_{13}+\lambda e_{24}$ takes the form
\begin{equation}\label{CondInv A_i:e_{13}+lambda e_{24}-3,19}
\begin{split}
&XA_1-A_3=0,\ \ XA_2+\lambda A_4=0,\\
&XA_3+A_1=0,\ \ XA_4+\lambda A_2=0,
\end{split}
\end{equation} %\eqno(3.129)
where
$X=x^3\partial_1+\lambda x^4\partial_2-x^1\partial_3+\lambda x^2\partial_4$.
By substitution \eqref{Change:P_{2,11a}} we replace $X$ to the form
\begin{equation}\label{X:e_{13}+lambda e_{24}-ch}
X=\lambda\tilde x^1\frac{\partial}{\partial\tilde x^1}-
 (\tilde x^3+\lambda\tilde x^2)\frac{\partial}{\partial\tilde x^2}+
 (\tilde x^2-\lambda\tilde x^3)\frac{\partial}{\partial\tilde x^3}.
\end{equation} %\eqno(3.130)
Substituting \eqref{Solution:P_{2,11a}} for $A_i$ in
\eqref{CondInv A_i:e_{13}+lambda e_{24}-3,19}--%
\eqref{X:e_{13}+lambda e_{24}-ch}, we get some equation;
using a linear independence of variables $\tilde x^2$, $\tilde x^3$, and
their powers, we obtain the system of equations:
\begin{equation}\label{Equations:P{3,19}}
\begin{split}
&\lambda\left(\tilde x^1\frac{\partial\Phi}{\partial\tilde x^1}-
 \Phi\right)=0,\ \ \ \
 \lambda\tilde x^1\frac{\partial\Psi}{\partial\tilde x^1}-\Xi=0,\\
&\lambda\left(\tilde x^1\frac{\partial\Theta}{\partial\tilde x^1}+
 \Theta-\Phi\right)=0,\ \
 \lambda\tilde x^1\frac{\partial\Xi}{\partial\tilde x^1}+\Psi=0.
\end{split}
\end{equation} %\eqno(3.131)
The solution of equations \eqref{Equations:P{3,19}} takes the form
\begin{equation}\label{Phi,Psi,Xi,Theta:P{3,19}}
\begin{split}
&\Phi=\tilde x^1C_1(\tilde x^4),\ \
 \Theta=\frac{\tilde x^1C_1(\tilde x^4)}2+\frac{C_2(\tilde x^4)}{\tilde x^1},\\
&\Psi=C_3(\tilde x^4)\cos\frac{\ln\tilde x^1}{\lambda}+
 C_4(\tilde x^4)\sin\frac{\ln\tilde x^1}{\lambda},\\
&\Xi=-C_3(\tilde x^4)\sin\frac{\ln\tilde x^1}{\lambda}+
 C_4(\tilde x^4)\cos\frac{\ln\tilde x^1}{\lambda},\\
\end{split}
\end{equation} %\eqno(3.132)
where $C_k(\tilde x^4)$ are arbitrary functions.
\begin{state} \label{state:P_{3,19}}
The class $P_{3,19}$ of potentials that admit the group $G_{3,19}$
is defined by formulae \eqref{Solution:P_{2,11a}} and
\eqref{Phi,Psi,Xi,Theta:P{3,19}}
(the transformation of coordinates is defined by \eqref{Change:P_{2,11a}}).
\end{state}

\subsubsection{Class $P_{3,20}$}  % subsubsection 3.3.20
The algebra $\mathcal L_{3,20}=L\{e_{12},\,e_{13},\,e_{23}\}$ corresponds
to the group $G_{3,20}=SO(3)$ of rotations over origin $O$ in the
three-di\-men\-si\-onal subspace $\mathbb R^3_0=\{x\in\mathbb R^4_1:\,x^4=0\}$
of Minkowski space. Since for $\lambda=\mu=0$
$\mathcal L_{1,2}\subset \mathcal L_{3,20}$ and the algebra
$\mathcal L_{3,20}$ is an extension of the algebra $L\{e_{13}\}$ by means of
the vectors $e_{12}$ and $e_{23}$, hence $P_{3,20}\subset P_{1,2}$
(for $\lambda=\mu=0$). The equation \eqref{CondInv A_i} for the vectors
$e_{12}$ and $e_{23}$ takes the forms
\begin{align}\label{CondInv A_i:e_{12}}
XA_1+A_2=0,\ \ XA_2-A_1=0,\ \ XA_3=0,\ \ XA_4=0
\end{align} %\eqno(3.133)
$(X=-x^2\partial_1+x^1\partial_2)$ and
\begin{align}\label{CondInv A_i:e_{23}}
YA_1=0,\ \ YA_2+A_3=0,\ \ YA_3-A_2=0,\ \ YA_4=0
\end{align} %\eqno(3.134)
$(Y=-x^3\partial_2+x^2\partial_3)$. By substitution \eqref{Change:Elipt}
for $\lambda=\mu=0$
\begin{equation}\label{Change:rot-e_{13}-gl2}
x^1=r\sin\varphi,\ \ x^2=\tilde x^2,\ \
x^3=r\cos\varphi,\ \ x^4=\tilde x^4
\end{equation} %\eqno(3.135)
we replace $X$ and $Y$ to the forms:
\begin{align}\label{X:e_{12}-gl2-ch}
&X=-\tilde x^2\sin\varphi\frac{\partial}{\partial r}+
 r\sin\varphi\frac{\partial}{\partial\tilde x^2}-
 \frac{\tilde x^2\cos\varphi}r\frac{\partial}{\partial\varphi},\\%\eqno(3.136)
    \label{Y:e_{23}-gl2-ch}
&Y=\tilde x^2\cos\varphi\frac{\partial}{\partial r}-
 r\cos\varphi\frac{\partial}{\partial\tilde x^2}-
 \frac{\tilde x^2\sin\varphi}r\frac{\partial}{\partial \varphi}. %\eqno(3.137)
\end{align}
Substituting \eqref{Cl:Elipt} for $A_i$ in equations
\eqref{CondInv A_i:e_{12}}--\eqref{X:e_{12}-gl2-ch} and
\eqref{CondInv A_i:e_{23}}--\eqref{Y:e_{23}-gl2-ch}, we get some equations;
dividing variables in them, we obtain the system of equations:
\begin{equation}\label{Equations:P{3,20}}
\begin{split}
&-\tilde x^2\frac{\partial C_2}{\partial r}+
 r\frac{\partial C_2}{\partial\tilde x^2}+A_2=0,\ \
 -\tilde x^2\frac{\partial C_1}{\partial r}+
 r\frac{\partial C_1}{\partial\tilde x^2}+\frac{C_1\tilde x^2}r=0,\\
&-\frac{C_2\tilde x^2}r+A_2=0,\ \
 -\tilde x^2\frac{\partial A_2}{\partial r}+
 r\frac{\partial A_2}{\partial\tilde x^2}-C_2=0,\ \ C_1=0,\\
&-\tilde x^2\frac{\partial C_1}{\partial r}+
 r\frac{\partial C_2}{\partial\tilde x^2}+\frac{C_2\tilde x^2}r=0,\ \
 -\tilde x^2\frac{\partial A_4}{\partial r}+
 r\frac{\partial A_4}{\partial\tilde x^2}=0.
\end{split}
\end{equation} %\eqno(3.138)
The solution of the system \eqref{Equations:P{3,20}} takes the form
\begin{equation}\label{Cl:P{3,20}-C_k}
C_1=C_2=A_2=0,\ \ A_4=A_4(\rho,\,x^4),
\end{equation} %\eqno(3.139)
where
$\rho=\sqrt{r^2+(\tilde x^2)^2}=\sqrt{(x^1)^2+(x^2)^2+(x^3)^2}$.
Substituting \eqref{Cl:P{3,20}-C_k} for $C_k$ and $A_i$ in \eqref{Cl:Elipt},
we get the following result.
\begin{state} \label{state:P_{3,20}}
The class $P_{3,20}$ of potentials that admit the group $G_{3,20}$
takes the form
\begin{equation}\label{Cl:P{3,20}}
A_i=(0,\ 0,\ 0,\ A_4(\rho,\,x^4)),
\end{equation} %\eqno(3.140)
where $A_4(\rho,\,x^4)$ is an arbitrary function.
\end{state}

\subsubsection{Class $P_{3,21}$}  % subsubsection 3.3.21
The algebra $\mathcal{L}_{3,21}=L\{e_{12},\,e_{14},\,e_{24}\}$ corresponds
to the group $G_{3,21}$ generated by rotations and by pseudo-rotations.
The algebra $\mathcal{L}_{3,21}$ is an extension of the algebra
$\mathcal{L}_{1,3}$ ($\lambda=0$) by means of the vectors $e_{12}$ and
$e_{14}$, hence $P_{3,21}\subset P_{1,3}$ (for $\lambda=0$).
The equation \eqref{CondInv A_i} for the vector $e_{12}$ takes the form
\eqref{CondInv A_i:e_{12}}, and for the vector $e_{14}$ --- as follows
\begin{equation}\label{CondInv A_i:e_{14}}
x^4\partial_1A_i+x^1\partial_4A_i+A_1\delta^4_i+A_4\delta^1_i=0.
\end{equation} %\eqno(3.141)
By substitution \eqref{Change:Hyperb1,3-gl2} for $\lambda=0$
\begin{equation}\label{Change:Hyperb3,21-gl2}
x^1=\tilde x^1,\ \ x^2=r\ch\varphi,\ \
x^3=\tilde x^3,\ \ x^4=r\sh\varphi
\end{equation} %\eqno(3.142)
we replace the operator $X=-x^2\partial_1+x^1\partial_2$ to the form:
\begin{equation}\label{X:e_{12}-gl2-3,21}
X=-r\ch\varphi\frac{\partial}{\partial \tilde x^1}+
 \tilde x^1\ch\varphi\frac{\partial}{\partial r}-
\frac{\tilde x^1\sh\varphi}r\frac{\partial}{\partial \varphi}.
\end{equation} %\eqno(3.143)
Substituting \eqref{Cl:P{1,3}-gl2} for $A_i$ in equations
\eqref{CondInv A_i:e_{12}}--\eqref{X:e_{12}-gl2-3,21}, we get some equations;
dividing variables in them, we obtain the system of equations:
\begin{equation}\label{Equations:P{3,21}}
\begin{split}
&C_2=0,\ \ r\frac{\partial A_1}{\partial\tilde x^1}-
 \tilde x^1\frac{\partial A_1}{\partial r}-C_1=0,\\
&r\frac{\partial C_1}{\partial\tilde x^1}-
 \tilde x^1\frac{\partial C_1}{\partial r}+A_1=0,\ \
 A_1+\frac{\tilde x^1C_1}r=0,\\
&r\frac{\partial A_3}{\partial\tilde x^1}-
 \tilde x^1\frac{\partial A_3}{\partial r}=0,\ \
 r\frac{\partial C_1}{\partial\tilde x^1}-
 \tilde x^1\frac{\partial C_1}{\partial r}+\frac{\tilde x^1C_1}r=0.
\end{split}
\end{equation} %\eqno(3.144)
The solution of the system \eqref{Equations:P{3,21}} takes the form
\begin{equation}\label{Cl:P{3,21}}
A_1=A_2=A_4=0,\ \ A_3=A_3(u,\,x^3),
\end{equation} %\eqno(3.145)
where $u=\sqrt{r^2+(\tilde x^1)^2}=\sqrt{(x^1)^2+(x^2)^2-(x^4)^2}$.
Substituting \eqref{Cl:P{3,21}} for $A_i$ in \eqref{CondInv A_i:e_{14}},
we get an identity. Thus we obtain the following result.
\begin{state} \label{state:P_{3,21}}
The class $P_{3,21}$ of potentials that admit the group $G_{3,21}$
takes the form \eqref{Cl:P{3,21}}, where $A_3(u,\,x^3)$ is an arbitrary
function.
\end{state}

%%%-------------------------------------------------------------

\subsection{Potentials that admit four-dimensional symmetry groups}
%Subsection  3.4
%
\subsubsection{Class $P_{4,1}$}  % subsubsection 3.4.1
The algebra $\mathcal{L}_{4,1}=L\{e_1,e_2,e_3,e_4\}$ corresponds to the group
of translations of Minkowski space $\mathbf{R}^4_1$. The algebra
$\mathcal{L}_{4,1}$ is an extension of the algebra $\mathcal{L}_{3,1a}$
by means of the vector $e_4$, hence $P_{4,1}\subset P_{3,1a}$.
Substituting $A_i(x^4)$ for $A_i$ in \eqref{CondInv A_i:e_{14}}, we get
\begin{state} \label{Cl:P{4,1}}
The class $P_{4,1}$ of potentials that admit the group $G_{4,1}$
consists of the fields $A_i$, which are constant in Galilean coordinates:
$A_i=\text{\rm const}$.
\end{state}

\subsubsection{Class $P_{4,2}$}  % subsubsection 3.4.2
The algebra $\mathcal{L}_{4,2}=L\{e_{13}+\mu e_4,\,e_1,\,e_2,\,e_3\}$
($\mu\ne0$) is an extension of the algebra $\mathcal{L}_{3,3}$
by means of the vector $e_2$, therefore $P_{4,2}\subset P_{3,3}$.
Substituting \eqref{Cl:P_{3,3}} for $A_i$ in \eqref{CondInv A_i:e_2}, we get
\begin{equation}\label{Cl:P_{4,2}}
\begin{split}
&A_1=C_1\sin\frac{x^4}{\mu}+C_2\cos\frac{x^4}{\mu},\ \ A_2=C_3,\\
&A_3=C_1\cos\frac{x^4}{\mu}-C_2\sin\frac{x^4}{\mu},\ \ A_4=C_4\ \
(C_k=\text{\rm const}).
\end{split}
\end{equation} %\eqno(3.146)
\begin{state} \label{Cl:P{4,2}}
The class $P_{4,2}$ of potentials that admit the group $G_{4,2}$
consists of the fields \eqref{Cl:P_{4,2}}.
\end{state}

\subsubsection{Class $P_{4,3}$}  % subsubsection 3.4.3
The algebra $\mathcal{L}_{4,3}=L\{e_{13}+\lambda e_2,\,e_1,\,e_3,\,e_4\}$
($\lambda\ne0$) is an extension of the algebra $\mathcal{L}_{3,2}$
by means of the vector $e_4$, therefore $P_{4,3}\subset P_{3,2}$. Substituting
\eqref{Cl:P_{3,2}-ne0} for $A_i$ in \eqref{CondInv A_i:e_4}, we get
\begin{equation}\label{Cl:P_{4,3}-ne0}
\begin{split}
&A_1=C_1\sin\frac{x^2}{\lambda}+C_2\cos\frac{x^2}{\lambda},\ \  A_2=C_3,\\
&A_3=C_1\cos\frac{x^2}{\lambda}-C_2\sin\frac{x^2}{\lambda},\ \  A_4=C_4,
\end{split}
\end{equation} %\eqno(3.147)
where $C_k=\text{\rm const}$. For $\lambda=0$ potentials of the class
$P_{3,2}$ are defined by \eqref{Cl:P_{3,2}-0}. Substituting
\eqref{Cl:P_{3,2}-0} for $A_i$ in \eqref{CondInv A_i:e_4}, we get
\begin{equation}\label{Cl:P_{4,3}-0}
A_1=A_3=0,\ \ A_2=A_2(x^2),\ \ A_4=A_4(x^2).
\end{equation} %\eqno(3.148)
\begin{state} \label{state:P{4,3}}
For $\lambda\ne0$ the class $P_{4,3}$ of potentials that admit the group
$G_{4,3}$ consists of the fields \eqref{Cl:P_{4,3}-ne0}; for $\lambda=0$ it
consists of the fields \eqref{Cl:P_{4,3}-0}.
\end{state}

\subsubsection{Class $P_{4,4}$}  % subsubsection 3.4.4
The algebra
$\mathcal{L}_{4,4}=L\{e_{13}+\lambda e_2,\ e_1,\,e_3,\,e_2+e_4\}$
is an extension of the algebra $\mathcal{L}_{3,1c}$ by means of the vector
$e_{13}+\lambda e_2$, therefore $P_{4,4}\subset P_{3,1c}$. Substituting
$A_i(x^2-x^4)$ for $A_i$ in \eqref{CondInv A_i-el(mu=0)}
(\eqref{CondInv A_i} for the vector $\xi=e_{13}+\lambda e_2$), we obtain
the system \eqref{P_{3,2}-sist}; for $\lambda\ne0$ the solution
of \eqref{P_{3,2}-sist} takes the form
\begin{equation}\label{Cl:P_{4,4}-ne0}
\begin{split}
&A_1=C_1\sin\frac{x^2-x^4}{\lambda}+C_2\cos\frac{x^2-x^4}{\lambda},\ \
  A_2=C_3,\\
&A_3=C_1\cos\frac{x^2-x^4}{\lambda}-C_2\sin\frac{x^2-x^4}{\lambda},\ \
  A_4=C_4,
\end{split}
\end{equation} %\eqno(3.149)
where $C_k=\text{\rm const}$. For $\lambda=0$ and $A_i=A_i(x^2-x^4)$
the solution of \eqref{P_{3,2}-sist} takes the form
\begin{equation}\label{Cl:P_{4,4}-0}
A_1=A_3=0,\ \ A_2=A_2(x^2-x^4),\ \ A_4=A_4(x^2-x^4),
\end{equation} %\eqno(3.150)
where $A_2(x^2-x^4)$ and $A_4(x^2-x^4)$ are arbitrary functions.
\begin{state} \label{state:P{4,4}}
For $\lambda\ne0$ the class $P_{4,4}$ of potentials that admit the group
$G_{4,4}$ consists of the fields \eqref{Cl:P_{4,4}-ne0}; for $\lambda=0$ it
consists of the fields \eqref{Cl:P_{4,4}-0}.
\end{state}

\subsubsection{Class $P_{4,5}$}  % subsubsection 3.4.5
The algebra $\mathcal{L}_{4,5}=L\{e_{24},\,e_1,\,e_3,\,e_2+e_4\}$
is an extension of the algebra $\mathcal{L}_{3,1c}$ by means of the vector
$e_{24}$, hence $P_{4,5}\subset P_{3,1c}$. Substituting $A_i(x^2-x^4)$
for $A_i$ in \eqref{CondInv A_i:e_{24}} (\eqref{CondInv A_i} for the vector
$\xi=e_{24}$), we obtain
\begin{equation}\label{P_{4,5}-sist}
\begin{split}
&(x^4-x^2)A'_1=0,\ \ (x^4-x^2)A'_2+A_4=0,\\
&(x^4-x^2)A'_3=0,\ \ (x^4-x^2)A'_4+A_2=0.\ \
\end{split}
\end{equation} %\eqno(3.151)
The solution of \eqref{P_{4,5}-sist} takes the form
\begin{equation}\label{Cl:P_{4,5}}
\begin{split}
&A_1=C_1,\ \ A_2=C_2\cdot(x^2-x^4)+\frac{C_4}{x^2-x^4},\\
&A_3=C_3,\ \ A_4=C_2\cdot(x^2-x^4)-\frac{C_4}{x^2-x^4},\ \
(C_k=\text{\rm const}).
\end{split}
\end{equation} %\eqno(3.152)
\begin{state} \label{state:P{4,5}}
The class $P_{4,5}$ of potentials that admit the group $G_{4,5}$
consists of the fields \eqref{Cl:P_{4,5}}.
\end{state}

\subsubsection{Class $P_{4,6}$}  % subsubsection 3.4.6
The algebra $\mathcal{L}_{4,6}=L\{e_{24}+\lambda e_3,\,e_1,\,e_2,\,e_4\}$
is an extension of the algebra $\mathcal{L}_{3,6}$ by means of the vector
$e_1$, hence $P_{4,6}\subset P_{3,6}$. Substituting \eqref{Cl:P_{3,6}-ne0}
for $A_i$ in \eqref{CondInv A_i:e_1} (\eqref{CondInv A_i} for the vector
$\xi=e_1$), we obtain for ${\lambda\ne0}$ the following result:
\begin{equation}\label{Cl:P_{4,6}-ne0}
\begin{split}
&A_1=C_1,\ \ A_2=C_2\ch\frac{x^3}{\lambda}+C_4\sh\frac{x^3}{\lambda},\\
&A_3=C_3,\ \ A_4=-C_2\sh\frac{x^3}{\lambda}-C_4\ch\frac{x^3}{\lambda},
\end{split}
\end{equation} %\eqno(3.153)
where $C_k=\text{\rm const}$.
Let be $\lambda=0$. Substituting \eqref{Cl:P_{3,6}-0} for $A_i$ in
\eqref{CondInv A_i:e_1}, we get
\begin{equation}\label{Cl:P_{4,6}-0}
A_1=A_1(x^3),\ \ A_2=0,\ \ A_3=A_3(x^3),\ \ A_4=0.
\end{equation} %\eqno(3.154)
\begin{state} \label{state:P{4,6}}
For $\lambda\ne0$ the class $P_{4,6}$ of potentials that admit the group
$G_{4,6}$ consists of the fields \eqref{Cl:P_{4,6}-ne0}; for $\lambda=0$ it
consists of the fields \eqref{Cl:P_{4,6}-0}.
\end{state}

\subsubsection{Class $P_{4,7}$}  % subsubsection 3.4.7
The algebra
$\mathcal{L}_{4,7}=L\{e_{13}+\lambda e_{24},\,e_1,\,e_3,\,e_2+e_4\}$
(${\lambda\ne0}$) is an extension of the algebra $\mathcal{L}_{3,1c}$
by means of the vector $e_{13}+\lambda e_{24}$, hence
$P_{4,7}\subset P_{3,1c}$. Substituting $A_i(x^2-x^4)$ for $A_i$ in
\eqref{CondInv A_i:birot}--\eqref{X:birot-gl2}
(\eqref{CondInv A_i} for the vector $\xi=e_{13}+\lambda e_{24}$), we obtain
\begin{equation}\label{Equations:P_{4,7}}
\begin{split}
&\lambda(x^4-x^2)A'_1-A_3=0,\ \ \lambda(x^4-x^2)A'_2+\lambda A_4=0,\\
&\lambda(x^4-x^2)A'_3+A_1=0,\ \ \lambda(x^4-x^2)A'_4+\lambda A_2=0.
\end{split}
\end{equation} %\eqno(3.155)
The solution of \eqref{Equations:P_{4,7}} takes the form
\begin{equation}\label{Cl:P_{4,7}}
\begin{split}
&A_1=C_1\cos\frac{\ln(x^2-x^4)}{\lambda}+C_3\sin\frac{\ln(x^2-x^4)}{\lambda},\\
&A_2=C_2(x^2-x^4)+\frac{C_4}{x^2-x^4},\\
&A_3=C_1\sin\frac{\ln(x^2-x^4)}{\lambda}-C_3\cos\frac{\ln(x^2-x^4)}{\lambda},\\
&A_4=C_2(x^2-x^4)-\frac{C_4}{x^2-x^4}\ \ \ \ (C_k=\text{\rm const}).
\end{split}
\end{equation} %\eqno(3.156)
\begin{state} \label{state:P{4,7}}
The class $P_{4,7}$ of potentials that admit the group $G_{4,7}$
consists of the fields \eqref{Cl:P_{4,7}}.
\end{state}

\subsubsection{Class $P_{4,8}$}  % subsubsection 3.4.8
The algebra
$\mathcal{L}_{4,8}=L\{e_{12}-e_{14}+\lambda e_3,\,e_1,\,e_2,\,e_4\}$
is an extension of the algebra $\mathcal{L}_{3,1b}$ by means of the vector
$e_{12}-e_{14}+\lambda e_3$, hence $P_{4,8}\subset P_{3,1b}$.
Substituting $A_i(x^3)$ for $A_i$ in the equation \eqref{CondInv A_i}
for the vector $\xi=e_{12}-e_{14}+\lambda e_3$
\begin{align}\label{CondInv A_i:prb-4,8}
&XA_i-A_1(\delta^2_i+\delta^4_i)+(A_2-A_4)\delta^1_i=0,\\ %\eqno(3.157)
\label{X:prb-P_{4,8}}
&X=-(x^2+x^4)\partial_1+x^1\partial_2+\lambda\partial_3-x^1\partial_4,
%\eqno(3.158)
\end{align}
we obtain
\begin{equation}\label{Equations:P_{4,8}}
\lambda A'_1+A_2-A_4=0,\ \lambda A'_2-A_1=0,\
\lambda A'_3=0,\ \lambda A'_4-A_1=0.
\end{equation} %\eqno(3.159)
For $\lambda\ne0$ the solution of \eqref{Equations:P_{4,8}} takes the form
\begin{equation}\label{Cl:P_{4,8}-ne0}
\begin{split}
&A_1=\frac{C_2}{\lambda}x^3+C_3,\ \
 A_2=\frac{C_2}{2\lambda^2}(x^3)^2+\frac{C_3}{\lambda}x^3+C_4,\\
&A_3=C_1,\ \ A_4=A_2+C_2\ \ \ \ (C_k=\text{\rm const});
\end{split}
\end{equation} %\eqno(3.160)
for $\lambda=0$ the same is
\begin{equation}\label{Cl:P_{4,8}-0}
A_1=0,\ \ A_2=A_4=\Phi(x^3),\ \ A_3=\Psi(x^3),
\end{equation} %\eqno(3.161)
where $\Phi(x^3)$ and $\Psi(x^3)$ are arbitrary functions.
\begin{state} \label{state:P{4,8}}
For $\lambda\ne0$ the class $P_{4,8}$ of potentials that admit the group
$G_{4,8}$ consists of the fields \eqref{Cl:P_{4,8}-ne0}; for $\lambda=0$ it
consists of the fields \eqref{Cl:P_{4,8}-0}.
\end{state}

\subsubsection{Class $P_{4,9}$}  % subsubsection 3.4.9
The algebra
$\mathcal{L}_{4,9}=$ $=L\{e_{12}-e_{14}+\lambda e_2,\,e_1,\,e_3,\,e_2-e_4\}$
corresponds to the group $G_{4,9}$ generated by parabolic helices and by
translations along the vectors of an isotropic hyperplane. The solution of
the system $L_{e_1}A_i=L_{e_3}A_i=L_{e_2-e_4}A_i=0$ takes the form
\begin{equation}\label{Cl:P_{3,1d}}
A_i=A_i(x^2+x^4).
\end{equation} %\eqno(3.162)
Substituting \eqref{Cl:P_{3,1d}} for $A_i$ in \eqref{CondInv A_i} for the
vector $e_{12}-e_{14}+\lambda e_2$
\begin{align}\label{CondInv A_i:prb{4,9}}
&XA_i-A_1(\delta^2_i+\delta^4_i)+(A_2-A_4)\delta^1_i=0,\\ %\eqno(3.163)
\label{X:prb-gl2{4,9}}
&X=-(x^2+x^4)\partial_1+(x^1+\lambda)\partial_2-x^1\partial_4, %\eqno(3.164)
\end{align}
we get the system \eqref{Equations:P_{4,8}}. For $\lambda\ne0$ and
$A_i=A_i(x^2+x^4)$ the solution of \eqref{Equations:P_{4,8}} takes the form
\begin{align}\label{Cl:P_{4,9}-ne0}
&A_1=\frac{C_2}{\lambda}(x^2+x^4)+C_3,\ \
 A_2=\frac{C_2}{2\lambda^2}(x^2+x^4)^2+
 \frac{C_3}{\lambda}(x^2+x^4)+C_4,\nonumber\\
&A_3=C_1,\ \ A_4=A_2+C_2\ \ (C_k=\text{\rm const}); %\eqno(3.165)
\end{align}
for $\lambda=0$ the same is
\begin{equation}\label{Cl:P_{4,9}-0}
A_1=0,\ \ A_2=A_4=\Phi(x^2+x^4),\ \ A_3=\Psi(x^2+x^4),
\end{equation} %\eqno(3.166)
where $\Phi(x^2+x^4)$ and $\Psi(x^2+x^4)$ are arbitrary functions.
\begin{state} \label{state:P{4,9}}
For $\lambda\ne0$ the class $P_{4,9}$ of potentials that admit the group
$G_{4,9}$ consists of the fields \eqref{Cl:P_{4,9}-ne0}; for $\lambda=0$ it
consists of the fields \eqref{Cl:P_{4,9}-0}.
\end{state}

\subsubsection{Class $P_{4,10}$}  % subsubsection 3.4.10
The algebra ${\mathcal{L}_{4,10}=L\{e_{13},\,e_{24},\,e_1,\,e_3,\}}$
is an extension of the algebra $\mathcal{L}_{3,5}$ by means of the vector
$e_{13}$, hence ${P_{4,10}\subset P_{3,5}}$. Substituting \eqref{Cl:P_{3,5}}
for $A_i$ in \eqref{CondInv A_i} for the vector $\xi=e_{13}$
\begin{equation}\label{CondInv A_i:e_{13}-4,10}
x^3\partial_1A_i-x^1\partial_3A_i+A_1\delta^3_i-A_3\delta^1_i=0,
\end{equation} %\eqno(3.167)
we get some system; using the substitution \eqref{Change:P_{3,5}},
we obtain the following result
\begin{equation}\label{Cl:P_{4,10}}
\begin{split}
&A_1=0,\ \ A_2=C_1(\rho)\ch\varphi+C_2(\rho)\sh\varphi,\\
&A_3=0,\ \ A_4=-C_1(\rho)\sh\varphi-C_2(\rho)\ch\varphi,
\end{split}
\end{equation} %\eqno(3.168)
where $C_1(\rho)$ and $C_2(\rho)$ are arbitrary functions.
\begin{state} \label{state:P{4,10}}
The class $P_{4,10}$ of potentials that admit the group $G_{4,10}$
is defined by \eqref{Cl:P_{4,10}} and \eqref{Change:P_{3,5}}.
\end{state}

\subsubsection{Class $P_{4,11}$}  % subsubsection 3.4.11
The algebra $\mathcal{L}_{4,11}=L\{e_{13},\,e_{24},\,e_2,\,e_4,\}$
is an extension of the algebra $\mathcal{L}_{3,6}$ (for $\lambda=0$)
by means of the vector $e_{13}$, hence ${P_{4,11}\subset P_{3,6}}$
($\lambda=0$). Substituting \eqref{Cl:P_{3,6}-0} for $A_i$ in
\eqref{CondInv A_i:e_{13}-4,10}, we get some system; using the substitution
\begin{equation}\label{Change:P_{4,11}}
x^1=r\sin\varphi,\ \ x^3=r\cos\varphi,
\end{equation} %\eqno(3.169)
we obtain the following result
\begin{equation}\label{Cl:P_{4,11}}
\begin{split}
&A_1=C_1(r)\sin\varphi+C_2(r)\cos\varphi,\ \ A_2=0,\\
&A_3=C_1(r)\cos\varphi-C_2(r)\sin\varphi,\ \ A_4=0,
\end{split}
\end{equation} %\eqno(3.170)
where $C_1(r)$ and $C_2(r)$ are arbitrary functions.
\begin{state} \label{state:P{4,11}}
The class $P_{4,11}$ of potentials that admit the group $G_{4,11}$
is defined by \eqref{Cl:P_{4,11}} and \eqref{Change:P_{4,11}}.
\end{state}

\subsubsection{Class $P_{4,12}$}  % subsubsection 3.4.12
The algebra
$$
\mathcal{L}_{4,12}=L\{e_{12}-e_{14}+\mu e_3,\ e_{23}+e_{34}+\nu e_2,\
e_1,\ e_2-e_4\}
$$
corresponds to the group $G_{4,12}$ generated by two one-di\-men\-si\-onal
subgroups of parabolic helices and by translations along the vectors of
an isotropic two-di\-men\-si\-onal plane. The total solution of the system of
equations \eqref{CondInv A_i} for $e_1$ and $e_2-e_4$ takes the form
\begin{equation}\label{Cl:P_{2,1d}}
A_i=A_i(x^2+x^4,\,x^3).
\end{equation} %\eqno(3.171)
The equation \eqref{CondInv A_i} for the vector $\xi=e_{12}-e_{14}+\mu e_3$
takes the form
\begin{gather}
XA_i-A_1(\delta^2_i+\delta^4_i)+(A_2-A_4)\delta^1_i=0,
\label{CondInv A_i:e_{12}-e_{14}+mu e_3}\\%\eqno(3.172)
X=-(x^2+x^4)\partial_1+x^1(\partial_2-\partial_4)+\mu\partial_3;
\label{X:e_{12}-e_{14}+mu e_3}%\eqno(3.173)
\end{gather}
for the vector $\xi=e_{23}+e_{34}+\nu e_2$ the same is
\begin{gather}
YA_i+A_3(\delta^2_i+\delta^4_i)-(A_2-A_4)\delta^3_i=0,
\label{CondInv A_i:e_{23}+e_{34}+nu e_2}\\%\eqno(3.174)
Y=\nu\partial_2+(x^2+x^4)\partial_3-x^3(\partial_2-\partial_4).
\label{Y:e_{23}+e_{34}+nu e_2}%\eqno(3.175)
\end{gather}
Substituting \eqref{Cl:P_{2,1d}} for $A_i$ in
\eqref{CondInv A_i:e_{12}-e_{14}+mu e_3}--\eqref{X:e_{12}-e_{14}+mu e_3} and
\eqref{CondInv A_i:e_{23}+e_{34}+nu e_2}--\eqref{Y:e_{23}+e_{34}+nu e_2},
we get
\begin{equation}\label{Equations:P_{4,12}-X}
\begin{split}
&\mu\frac{\partial A_1}{\partial x^3}+A_2-A_4=0,\ \
 \mu\frac{\partial A_2}{\partial x^3}-A_1=0,\\
&\mu\frac{\partial A_3}{\partial x^3}=0,\ \
 \mu\frac{\partial A_4}{\partial x^3}-A_1=0
\end{split}
\end{equation} %\eqno(3.176)
and
\begin{align}\label{Equations:P_{4,12}-Y}
&\nu\frac{\partial A_1}{\partial x^2}+
 (x^2+x^4)\frac{\partial A_1}{\partial x^3}=0,\ \
 \nu\frac{\partial A_2}{\partial x^2}+
 (x^2+x^4)\frac{\partial A_2}{\partial x^3}+A_3=0,\nonumber\\
&\nu\frac{\partial A_3}{\partial x^2}+
 (x^2+x^4)\frac{\partial A_3}{\partial x^3}-A_2+A_4=0,\nonumber\\
&\nu\frac{\partial A_4}{\partial x^2}+
 (x^2+x^4)\frac{\partial A_4}{\partial x^3}+A_3=0.
\end{align} %\eqno(3.177)
Solving the system of equations
\eqref{Equations:P_{4,12}-X}--\eqref{Equations:P_{4,12}-Y}, we obtain the
result.
\begin{state} \label{state:P{4,12}}
The class $P_{4,12}$ of potentials that admit the group $G_{4,12}$
is defined by the following formulae:
\newline
a) for $\mu=\nu=0$
\begin{equation}\label{Cl:P_{4,12}-a}
\begin{split}
&A_1=0,\ \ A_3=\Psi(x^2+x^4),\\
&A_2=A_4=\Phi(x^2+x^4)-\frac{x^3\Psi(x^2+x^4)}{x^2+x^4};
\end{split}
\end{equation} %\eqno(3.178)
b) for $\mu=0,\ \nu\ne0$
\begin{equation}\label{Cl:P_{4,12}-b}
\begin{split}
&A_1=0,\ \ A_2=A_4=\Phi(u)-\frac{x^2+x^4}{\nu}\Psi(u),\\
&A_3=\Psi(u)\ \ \left(u=x^3-\frac{(x^2+x^4)^2}{2\nu}\right);
\end{split}
\end{equation} %\eqno(3.179)
c) for $\mu\ne0,\ \nu=0$
\begin{equation}\label{Cl:P_{4,12}-c}
\begin{split}
&A_1=\Phi(x^2+x^4),\ \ A_3=-\frac{x^2+x^4}{\mu}\Phi(x^2+x^4),\\
&A_2=A_4=\frac{x^3}{\mu}\Phi(x^2+x^4)+\Psi(x^2+x^4);
\end{split}
\end{equation} %\eqno(3.180)
d) for $\mu\ne0,\ \nu\ne0$
\begin{equation}\label{Cl:P_{4,12}-d}
A_1=A_3=0,\ \ A_2=A_4=\text{\rm const}
\end{equation} %\eqno(3.181)
($\Phi$ and $\Psi$ are arbitrary functions of one variable).
\end{state}

\subsubsection{Class $P_{4,13}$}  % subsubsection 3.4.13
The algebra
$$
\mathcal{L}_{4,13}=L\{e_{12}-e_{14},\ e_{24}+\lambda e_1,\ e_3,\ e_2-e_4\}
$$
corresponds to the group $G_{4,13}$ generated by parabolic rotations, by
hyperbolic helices, and by translations along the vectors of an isotropic
two-di\-men\-si\-onal plane. The solution of the system of equations
\eqref{CondInv A_i} for $e_3$ and $e_2-e_4$ takes the form
\begin{equation}\label{Cl:P_{2,1e}}
A_i=A_i(x^1,\,x^2+x^4).
\end{equation} %\eqno(3.182)
Substituting \eqref{Cl:P_{2,1e}} for $A_i$ in equation \eqref{CondInv A_i}
for the vector $e_{12}-e_{14}$
\begin{gather}
XA_i-A_1(\delta^2_i+\delta^4_i)+(A_2-A_4)\delta^1_i=0,
\label{CondInv A_i:e_{12}-e_{14}}\\ %\eqno(3.183)
X=-(x^2+x^4)\partial_1+x^1(\partial_2-\partial_4),
\label{X:e_{12}-e_{14}} %\eqno(3.184)
\end{gather}
we have the result:
\begin{equation}\label{P:e_{12}-e_{14},e_3,e_2-e_4}
\begin{split}
&A_1=-\frac{x^1C_1}{x^2+x^4}+C_2,\ \
 A_2=\frac{(x^1)^2C_1}{2(x^2+x^4)}-\frac{x^1C_2}{x^2+x^4}+C_3,\\
&A_3=C_4,\ \ A_4=A_2+C_1\ \ (C_k=C_k(x^2+x^4)).
\end{split}
\end{equation} %\eqno(3.185)
Substituting \eqref{P:e_{12}-e_{14},e_3,e_2-e_4} for $A_i$ in
\eqref{CondInv A_i-Hyp} (the equation \eqref{CondInv A_i} for the vector
${e_{24}+\lambda e_1}$), we obtain
\begin{equation}\label{Cl:P_{4,13}-C_k}
\begin{split}
&C_1=K_1(x^2+x^4),\ \ C_2=K_1\lambda\ln(x^2+x^4)+K_2,\\
&C_3=\frac{K_1\lambda^2\ln^2(x^2+x^4)+2K_2\lambda\ln(x^2+x^4)}{2(x^2+x^4)}-\\
&-\frac{K_1}{2}(x^2+x^4)+\frac{K_3}{x^2+x^4},\ \ C_4=K_4\ \
(K_i=\text{\rm const)}.
\end{split}
\end{equation} %\eqno(3.186)
\begin{state} \label{state:P{4,13}}
The class $P_{4,13}$ of potentials that admit the group $G_{4,13}$
is defined by \eqref{P:e_{12}-e_{14},e_3,e_2-e_4} and \eqref{Cl:P_{4,13}-C_k}.
\end{state}

\subsubsection{}  % subsubsection 3.4.14 Class $P_{4,14}$
The algebra
$$
\mathcal{L}_{4,14}=L\{e_{12}-e_{14},\,e_{24}+\lambda e_3,\,
e_1+\nu e_3,\,e_2-e_4\}
$$
corresponds to the group $G_{4,14}$ generated by parabolic rotations, by
hyperbolic helices, and by translations along the vectors of an isotropic
two-di\-men\-si\-onal plane (the group $G_{4,14}$ is not conjugated
to $G_{4,13}$). The solution of the system of equations \eqref{CondInv A_i}
for $e_1+\nu e_3$ and $e_2-e_4$ takes the form
\begin{equation}\label{Cl:P_{2,1f}}
A_i=A_i(x^2+x^4,\,x^3-\nu x^1).
\end{equation} %\eqno(3.187)

3.4.14.1. {\it Class $P_{4,14a}$.} %% Point 3.4.14.1
Let be $\nu\ne0$. Substituting \eqref{Cl:P_{2,1f}} for $A_i$ in
\eqref{CondInv A_i:e_{12}-e_{14}}--\eqref{X:e_{12}-e_{14}} (the equation
\eqref{CondInv A_i} for the vector $e_{12}-e_{14}$), we obtain
\begin{equation}\label{P:e_{12}-e_{14},e_1+nu e_3,e_2-e_4}
\begin{split}
&A_1=\frac{(x^3-\nu x^1)C_1}{\nu(x^2+x^4)}+C_2,\\
&A_2=\frac{(x^3-\nu x^1)^2C_1}{2\nu^2(x^2+x^4)^2}+
 \frac{(x^3-\nu x^1)C_2}{\nu(x^2+x^4)}+C_3,\\
&A_3=C_4,\ \ A_4=A_2+C_1\ \ (C_k=C_k(x^2+x^4)).
\end{split}
\end{equation} %\eqno(3.188)
Substituting \eqref{P:e_{12}-e_{14},e_1+nu e_3,e_2-e_4} for $A_i$ in equation
\begin{equation}\label{CondInv A_i:e_{24}+nu e_3}
x^4\partial_2A_i+\nu\partial_3A_i+x^2\partial_4A_i+
A_2\delta^4_i+A_4\delta^2_i=0
\end{equation} %\eqno(3.189)
(\eqref{CondInv A_i} for the vector $\xi=e_{24}+\nu e_3$), we get
\begin{equation}\label{Cl:P_{4,14a}-C_k}
\begin{split}
&C_1=K_1(x^2+x^4),\ \ C_2=-\frac{K_1\lambda}{\nu}\ln(x^2+x^4)+K_2,\\
&C_3=\frac{K_1\lambda^2\ln^2(x^2+x^4)-2K_2\lambda\nu\ln(x^2+x^4)}
 {2\nu^2(x^2+x^4)}-\\
&-\frac{K_1}{2}(x^2+x^4)+\frac{K_3}{x^2+x^4},\ \ C_4=K_4\ \ (K_i=\text{const}).
\end{split}
\end{equation} %\eqno(3.190)
\begin{state} \label{state:P{4,14a}}
The class $P_{4,14a}$ of potentials that admit the group $G_{4,14a}$
is defined by \eqref{P:e_{12}-e_{14},e_1+nu e_3,e_2-e_4} and
\eqref{Cl:P_{4,14a}-C_k}.
\end{state}

3.4.14.2. {\it Class $P_{4,14b}$.} %% Point 3.4.14.2
For $\nu=0$ we have the following result.
\begin{state} \label{state:P{4,14b}}
The class $P_{4,14b}$ of potentials that admit the group $G_{4,14b}$
($G_{4,14}$ for $\nu=0$), is defined by the formulae
\begin{equation}\label{Cl:P_{4,14b}}
\begin{split}
&A_1=0,\ \ A_3=\Psi(x^3-\lambda\ln(x^2+x^4)),\\
&A_2=A_4=\frac{\Phi(x^3-\lambda\ln(x^2+x^4))}{x^2+x^4},
\end{split}
\end{equation} %\eqno(3.191)
where $\Phi$ and $\Psi$ are arbitrary functions of one variable.
\end{state}

\subsubsection{Class $P_{4,15}$}  % subsubsection 3.4.15
The algebra $\mathcal{L}_{4,15}=L\{e_{12}-e_{14},\ e_{23}+e_{34},\
{e_{24}+\lambda e_1},\break {e_2-e_4}\}$
is an extension of the algebra $\mathcal{L}_{2,11a}$ by means of the vectors
$e_{24}+\lambda e_1$ and $e_2-e_4$, hence $P_{4,15}\subset P_{2,11a}$.
Substituting \eqref{Solution:P_{2,11a}}--\eqref{Change:P_{2,11a}} for $A_i$
in equation
\begin{equation}\label{CondInv A_i:e_2-e_4-4,15}
L_{e_2-e_4}A_i=\partial_2A_i-\partial_4A_i=
2\tilde x^1\frac{\partial A_i}{\partial\tilde x^4}=0,
\end{equation} %\eqno(3.192)
we get
\begin{align}\label{Cl:P_{4,15}-A_i}
&A_1=\frac{x^1C_1}{x^2+x^4}+C_2,\ \
A_2=-\frac{((x^1)^2+(x^3)^2)C_1}{2(x^2+x^4)^2}-
 \frac{x^1C_2+x^3C_3}{x^2+x^4}+C_4,\nonumber\\
&A_3=\frac{x^3C_1}{x^2+x^4}+C_3,\ \ A_4=A_2-C_1\ \
 (C_k=C_k(x^2+x^4)).
\end{align} %\eqno(3.193)
Substituting \eqref{Cl:P_{4,15}-A_i} for $A_i$ in \eqref{CondInv A_i-Hyp}
(equation \eqref{CondInv A_i} for ${e_{24}+\lambda e_1}$), we obtain
\begin{equation}\label{Cl:P_{4,15}-C_k}
\begin{split}
&C_1=K_1(x^2+x^4),\ \ C_2=-K_1\lambda\ln(x^2+x^4)+K_2,\\
&C_3=K_3,\ \
C_4=\frac{-K_1\lambda^2\ln^2(x^2+x^4)+2K_2\lambda\ln(x^2+x^4)}{2(x^2+x^4)}+\\
&+\frac{K_1}{2}(x^2+x^4)+\frac{K_4}{x^2+x^4},\ \ (K_i=\text{\rm const}).
\end{split}
\end{equation} %\eqno(3.194)
\begin{state} \label{state:P{4,15}}
The class $P_{4,15}$ of potentials that admit the group $G_{4,15}$
is defined by \eqref{Cl:P_{4,15}-A_i} and \eqref{Cl:P_{4,15}-C_k}.
\end{state}

\subsubsection{Class $P_{4,16}$}  % subsubsection 3.4.16
The algebra
$$
\mathcal L_{4,16}=L\{e_{12}-e_{14}+\lambda e_3,\
e_{23}+e_{34}+\lambda e_1,\ e_{13},\ e_2-e_4\}
$$
is an extension of the algebra $\mathcal L_{3,14}$\footnotemark[1]%
\footnotetext[1]{If we replace ${\mu\mapsto\lambda}$, ${\nu\mapsto\lambda}$,
${\lambda\mapsto0}$.}
by means of the vector $e_{13}$; therefore $P_{4,16}\subset P_{3,14}$ for
corresponding values of parametres. The substitution \eqref{Change:P{3,14}}
is replaced by
\begin{equation}\label{Change:P{4,16}}
u=x^2+x^4,\ \ \varphi=\frac{\lambda x^1+ux^3}{u^2+\lambda^2},\ \
\psi=\frac{\lambda x^3-x^1u}{u^2+\lambda^2}.
\end{equation} %\eqno(3.195)
Substituting \eqref{Cl:P{3,14}}---\eqref{Change:P{4,16}} for $A_i$ in
\eqref{CondInv A_i:e_{13}-4,10} (equation \eqref{CondInv A_i} for the vector
$e_{13}$), we obtain the following result.
\begin{state} \label{state:P{4,16}}
The class $P_{4,16}$ of potentials that admit the group $G_{4,16}$
is defined by the formulae
\begin{equation}\label{Cl:P{4,16}}
A_1=A_3=0,\ \ A_2=A_4=\Phi(x^2+x^4),
\end{equation} %\eqno(3.196)
where $\Phi(u)$ is an arbitrary function.
\end{state}

\subsubsection{Class $P_{4,17}$}  % subsubsection 3.4.17
The algebra
$$
\mathcal{L}_{4,17}=L\{e_{12}-e_{14},\,e_{23}+e_{34},\,
e_{13}+\lambda e_{24},\,e_2-e_4\}
$$
is an extension of the algebra $\mathcal{L}_{3,19}$ by means of the vector
$e_2-e_4$; therefore $P_{4,17}\subset P_{3,19}$. Substituting
\eqref{Solution:P_{2,11a}}--\eqref{Phi,Psi,Xi,Theta:P{3,19}} for $A_i$ in
\eqref{CondInv A_i:e_2-e_4-4,15}, we get result
\begin{equation}\label{Cl:P_{4,17}-A_i}
\begin{split}
&A_1=K_1x^1+\Psi,\ \  A_2=-\frac{K_1((x^1)^2+(x^3)^2)}{2(x^2+x^4)}-\\
&-\frac{x^1\Psi+x^3\Xi}{x^2+x^4}+\frac{K_1}{2}(x^2+x^4)+\frac{K_2}{x^2+x^4},\\
&A_3=K_1x^3+\Xi,\ \ A_4=A_2-K_1(x^2+x^4),
\end{split}
\end{equation} %\eqno(3.197)
where
\begin{equation}\label{Cl:P_{4,17}-Psi,Xi}
\begin{split}
&\Psi=K_3\cos\frac{\ln(x^2+x^4)}{\lambda}+
      K_4\sin\frac{\ln(x^2+x^4)}{\lambda},\\
&\Xi=-K_3\sin\frac{\ln(x^2+x^4)}{\lambda}+
      K_4\cos\frac{\ln(x^2+x^4)}{\lambda}\ \ \ (K_i=\text{\rm const}).
\end{split}
\end{equation} %\eqno(3.198)
\begin{state} \label{state:P{4,17}}
The class $P_{4,17}$ of potentials that admit the group $G_{4,17}$
is defined by the formulae \eqref{Cl:P_{4,17}-A_i} and
\eqref{Cl:P_{4,17}-Psi,Xi}.
\end{state}

\subsubsection{Class $P_{4,18}$}  % subsubsection 3.4.18
The algebra $\mathcal{L}_{4,18}=L\{e_{12}, e_{13}, e_{23}, e_4\}$
is an extension of the algebra $\mathcal{L}_{3,20}$ by means of the vector
$e_4$, hence $P_{4,18}\subset P_{3,20}$. Substituting \eqref{Cl:P{3,20}}
for $A_i$ in \eqref{CondInv A_i:e_4}, we get the following result.
\begin{state} \label{state:P{4,18}}
The class $P_{4,18}$ of potentials that admit the group $G_{4,18}$
takes the form
\begin{equation}\label{Cl:P{4,18}}
A_i=(0,\ 0,\ 0,\ A_4(\rho)),
\end{equation} %\eqno(3.199)
where $\rho=\sqrt{(x^1)^2+(x^2)^2+(x^3)^2}$, and $A_4(\rho)$
is an arbitrary function.
\end{state}

\subsubsection{Class $P_{4,19}$}  % subsubsection 3.4.19
The algebra $\mathcal{L}_{4,19}=L\{e_{12},\ e_{14},\ e_{24},\ e_3\}$
is an extension of the algebra $\mathcal{L}_{3,21}$ by means of the vector
$e_3$, hence ${P_{4,19}\subset P_{3,21}}$. Substituting \eqref{Cl:P{3,21}}
for $A_i$ in \eqref{CondInv A_i:e_3}, we get the following result.
\begin{state} \label{state:P{4,19}}
The class $P_{4,19}$ of potentials that admit the group $G_{4,19}$
takes the form
\begin{equation}\label{Cl:P{4,19}}
A_i=(0,\,0,\,C(u),\,0),
\end{equation} %\eqno(3.200)
where $u=\sqrt{(x^1)^2+(x^2)^2-(x^4)^2}$, and $C(u)$ is an arbitrary function.
\end{state}

\subsubsection{Class $P_{4,20}$}  % subsubsection 3.4.20
The algebra
$\mathcal{L}_{4,20}=L\{e_{12}-e_{14},\,e_{23}+e_{34},\,e_{13},\,e_{24}\}$
is an extension of the algebra $\mathcal{L}_{3,17}$ by means of the vector
$e_{13}$, hence $P_{4,20}\subset P_{3,17}$. Substituting
\eqref{Solution:P_{2,11a}}--\eqref{Phi,Psi,Xi,Theta:P_{3,17}} for $A_i$ in
\eqref{CondInv A_i:e_{13}+lambda(e_2-e_4)-3,18}--%
\eqref{X:e_{13}+lambda(e_2-e_4)-ch} (for $\lambda=0$), we get a solution;
returning to coordinates $\{x^i\}$, we obtain the result.
\begin{state} \label{state:P{4,20}}
The class $P_{4,20}$ of potentials that admit the group $G_{4,20}$
takes the form
\begin{equation}\label{Cl:P{4,20}}
\begin{split}
&A_1=x^1C(\tilde x^4),\ \ A_3=x^3C(\tilde x^4),\\
&A_2=\frac{C(\tilde x^4)}2\left(x^2+x^4- \frac{(x^1)^2+(x^3)^2}{x^2+x^4}\right)+
 \frac{D(\tilde x^4)}{x^2+x^4},\\
&A_4=A_2-(x^2+x^4)C(\tilde x^4),
\end{split}
\end{equation} %\eqno(3.201)
where $\tilde x^4=(x^1)^2+(x^2)^2+(x^3)^2-(x^4)^2$, and $C(\tilde x^4)$,
$D(\tilde x^4)$ are arbitrary functions.
\end{state}

%%%-------------------------------------------------------------

\newpage
\subsection{Potentials that admit five-dimensional symmetry groups}
%Subsection  3.5
%
\subsubsection{Class $P_{5,1}$}  % subsubsection 3.5.1
The algebra $\mathcal{L}_{5,1}=L\{e_{24},\,e_1,\,e_2,\,e_3,\,e_4\}$
is an extension of the algebra $\mathcal{L}_{4,1}$ by means of the vector
$e_{24}$, hence ${P_{5,1}\subset P_{4,1}}$. Substituting
$A_i=\text{\rm const}$ for $A_i$ in \eqref{CondInv A_i-Hyp} (for $\lambda=0$),
we obtain the result.
\begin{state} \label{state:P{5,1}}
The class $P_{5,1}$ of potentials that admit the group $G_{5,1}$
takes the form
\begin{equation}\label{Cl:P{5,1}}
A_i=(0,\,A,\,0,\,B)\ \ (A,\,B=\text{\rm const}).
\end{equation} %\eqno(3.202)
\end{state}

\subsubsection{Class $P_{5,2}$}  % subsubsection 3.5.2
The algebra
$\mathcal{L}_{5,2}=L\{e_{13}+\lambda e_{24},\,e_1,\,e_2,\,e_3,\,e_4\}$
is an extension of the algebra $\mathcal{L}_{4,1}$ by means of the vector
$e_{13}+\lambda e_{24}$, hence $P_{5,2}\subset P_{4,1}$. Substituting
$A_i=\text{\rm const}$ for $A_i$ in
\eqref{CondInv A_i:birot}--\eqref{X:birot-gl2}
($\lambda\ne0$), we obtain $A_i=0$, i.~e. {\it the class $P_{5,2}$ is empty}.

\subsubsection{Class $P_{5,3}$}  % subsubsection 3.5.3
The algebra $\mathcal{L}_{5,3}=L\{e_{12}-e_{14},\,e_1,\,e_2,\,e_3,\,e_4\}$
is an extension of the algebra $\mathcal{L}_{4,1}$ by means of the vector
$e_{12}-e_{14}$, hence $P_{5,3}\subset P_{4,1}$. Substituting
$A_i=\text{\rm const}$ for $A_i$ in \eqref{CondInv A_i:prb}--\eqref{X:prb-gl2}
(${\lambda=\mu=0}$), we obtain the result.
\begin{state} \label{state:P{5,3}}
The class $P_{5,3}$ of potentials that admit the group $G_{5,3}$
takes the form
\begin{equation}\label{Cl:P{5,3}}
A_i=(0,\,A,\,B,\,A)\ \ (A,\,B=\text{\rm const}).
\end{equation} %\eqno(3.203)
\end{state}

\subsubsection{Class $P_{5,4}$}  % subsubsection 3.5.4
The algebra $\mathcal{L}_{5,4}=L\{e_{13},\,e_{24},\,e_1,\,e_3,\,e_2+e_4\}$
is an extension of the algebra $\mathcal{L}_{4,5}$ by means of the vector
$e_{13}$, hence ${P_{5,4}\subset P_{4,5}}$. Substituting \eqref{Cl:P_{4,5}}
for $A_i$ in \eqref{CondInv A_i-el} ($\lambda=\mu=0$), we obtain
\begin{equation}\label{Cl:P_{5,4}}
\begin{split}
&A_1=0,\ \ A_2=K_1(x^2-x^4)+\frac{K_2}{x^2-x^4},\\
&A_3=0,\ \ A_4=K_1(x^2-x^4)-\frac{K_2}{x^2-x^4},
\end{split}
\end{equation} %\eqno(3.204)
where $K_1$ and $K_2$ are arbitrary constants.
\begin{state} \label{state:P{5,4}}
The class $P_{5,4}$ of potentials that admit the group $G_{5,4}$
is defined by \eqref{Cl:P_{5,4}}.
\end{state}

\subsubsection{Class $P_{5,5}$}  % subsubsection 3.5.5
The algebra
$$
\mathcal{L}_{5,5}=L\{e_{12}-e_{14},\ e_{23}+e_{34}+\lambda e_2,\
e_1,\ e_3,\ e_2-e_4\}
$$
is an extension of the algebra $\mathcal{L}_{4,9}$ ($\lambda=0$)
by means of the vector $e_{23}+e_{34}+\lambda e_2$, hence
$P_{5,4}\subset P_{4,9}$. Substituting \eqref{Cl:P_{4,9}-0} for $A_i$ in
\eqref{CondInv A_i:e_{23}+e_{34}+nu e_2}--\eqref{Y:e_{23}+e_{34}+nu e_2}
($\nu\mapsto\lambda$), we obtain the following result.
\begin{state} \label{state:P{5,5}}
For $\lambda=0$ the class $P_{5,5}$ of potentials that admit the group
$G_{5,5}$ takes the form
\begin{equation}\label{Cl:P{5,5}-0}
A_i=(0,\,\Phi(x^2+x^4),\,0,\,\Phi(x^2+x^4))
\end{equation} %\eqno(3.205)
where $\Phi$ is an arbitrary function; for $\lambda\ne0$ the same is
\begin{equation}\label{Cl:P_{5,5}-ne0}
A_1=0,\ \ A_2=A_4=-\frac{K_1}{\lambda}(x^2+x^4)+K_2,\ \ A_3=K_1,
\end{equation} %\eqno(3.206)
where $K_1$ and $K_2$ are arbitrary constants.
\end{state}

\subsubsection{Class $P_{5,6}$}  % subsubsection 3.5.6
The algebra
$
{\mathcal{L}_{5,6}=L\{e_{12}-e_{14},\,e_{24}+\lambda e_3,\,e_1,\,e_2,\,e_4\}}
$
is an extension of the algebra $\mathcal{L}_{4,6}$ by means of the vector
$e_{12}-e_{14}$, hence $P_{5,6}\subset P_{4,6}$. Substituting
\eqref{Cl:P_{4,6}-ne0} or \eqref{Cl:P_{4,6}-0} for $A_i$ in
\eqref{CondInv A_i:prb}--\eqref{X:prb-gl2} ($\lambda=\mu=0$), we get
\begin{equation}\label{Cl:P_{5,6}-ne0}
A_1=0,\ \ A_2=A_4=K_1e^{-x^3/\lambda},\ \ A_3=K_2\ \
(K_1,\ K_2=\text{\rm const})
\end{equation} %\eqno(3.207)
or
\begin{equation}\label{Cl:P_{5,6}-0}
A_i=(0,\,0,\,\Phi(x^3),\,0).
\end{equation} %\eqno(3.208)
\begin{state} \label{state:P{5,6}}
For $\lambda\ne0$ the class $P_{5,6}$ of potentials that admit the group
$G_{5,6}$ is defined by \eqref{Cl:P_{5,6}-ne0}; for $\lambda=0$ the same
is \eqref{Cl:P_{5,6}-0}.
\end{state}

\subsubsection{Class $P_{5,7}$}  % subsubsection 3.5.7
The algebra
$\mathcal{L}_{5,7}=L\{e_{12}-e_{14},\,e_{24},\,e_1,\,e_3,\,e_2-e_4\}$
is an extension of the algebra $\mathcal{L}_{4,13}$ by means of the vector
$e_1$, hence $P_{5,7}\subset P_{4,13}$. Substituting
\eqref{P:e_{12}-e_{14},e_3,e_2-e_4}--\eqref{Cl:P_{4,13}-C_k} for $A_i$ in
\eqref{CondInv A_i:e_1}, we obtain the following result.
\begin{state} \label{state:P{5,7}}
The class $P_{5,7}$ of potentials that admit the group $G_{5,7}$
takes the form
\begin{equation}\label{Cl:P_{5,7}}
A_i=\left(0,\,\frac{B}{x^2+x^4},\,C,\,\frac{B}{x^2+x^4}\right),
\end{equation} %\eqno(3.209)
where $B$ and $C$ are arbitrary constants.
\end{state}

\subsubsection{Class $P_{5,8}$}  % subsubsection 3.5.8
The algebra
$$
\mathcal{L}_{5,8}=L\{e_{12}-e_{14},\,e_{23}+e_{34},\,e_{24}+\lambda e_3,\,
e_1,\,e_2-e_4\}
$$
is an extension of the algebra $\mathcal{L}_{4,12a}$ ($\mathcal{L}_{4,12}$
for ${\mu=\nu=0}$) by means of the vector $e_{24}+\lambda e_3$, hence
$P_{5,8}\subset P_{4,12a}$. Substituting \eqref{Cl:P_{4,12}-a} for $A_i$ in
\eqref{CondInv A_i:e_{24}+lambda e_3}, we obtain the following result.
\begin{state} \label{state:P{5,8}}
The class $P_{5,8}$ of potentials that admit the group $G_{5,8}$
takes the form
\begin{equation}\label{Cl:P_{5,8}}
A_i=(0,\,\Phi(x^2+x^4),\,0,\,\Phi(x^2+x^4)),
\end{equation} %\eqno(3.210)
where $\Phi$ is an arbitrary function of one variable.
\end{state}

\subsubsection{Class $P_{5,9}$}  % subsubsection 3.5.9
The algebra
$$
\mathcal{L}_{5,9}=L\{e_{12}-e_{14},\ e_{23}+e_{34},\ e_{13},\
e_{24},\ e_2-e_4\}
$$
is an extension of the algebra $\mathcal{L}_{4,20}$ by means of the vector
$e_2-e_4$, hence $P_{5,9}\subset P_{4,20}$. Substituting
\eqref{Cl:P{4,20}} for $A_i$ in \eqref{CondInv A_i:e_2-e_4}, we get
\begin{equation}\label{Cl:P{5,9}}
\begin{split}
&A_1=Cx^1,\ \ A_2=\frac{C}2\left(x^2+x^4-
 \frac{(x^1)^2+(x^3)^2}{x^2+x^4}\right)+\frac{D}{x^2+x^4},\\
&A_3=Cx^3,\ \ A_4=A_2-C(x^2+x^4)\ \ (C,\,D=\text{\rm const}).
\end{split}
\end{equation} %\eqno(3.211)
\begin{state} \label{state:P{5,9}}
The class $P_{5,9}$ of potentials that admit the group $G_{5,9}$
is defined by \eqref{Cl:P{5,9}}.
\end{state}

%%%-------------------------------------------------------------

\subsection{Potentials that admit six-dimensional symmetry groups}
%Subsection  3.6
%
\subsubsection{Class $P_{6,1}$}  % subsubsection 3.6.1
The algebra
$\mathcal L_{6,1}=L\{e_{12},\,e_{13},\,e_{23},\,e_{14},\,e_{24},\,e_{34}\}$
corresponds to the Lorentz group. It is an extension of the algebra
$\mathcal{L}_{3,20}$ by means of the vectors $e_{14}$, $e_{24}$, and $e_{34}$,
hence $P_{6,1}\subset P_{3,20}$. Substituting \eqref{Cl:P{3,20}} for $A_i$ in
\eqref{CondInv A_i-Hyp} ($\lambda=0$), we obtain $A_i=0$, i.~e.
{\it the class of potentials that admit the Lorentz group $G_{6,1}$ is empty}.

\subsubsection{Class $P_{6,2}$}  % subsubsection 3.6.2
The algebra
${\mathcal{L}_{6,2}=L\{e_{13},\,e_{24},\,e_1,\,e_2,\,e_3,\,e_4\}}$
is an extension of the algebra $\mathcal{L}_{5,1}$ by means of the vector
$e_{13}$, hence $P_{6,2}\subset P_{5,1}$. Substituting \eqref{Cl:P{5,1}}
for $A_i$ in \eqref{CondInv A_i-el} ($\lambda=\mu=0$), we obtain $A_i=0$,
i.~e. {\it the class of potentials that admit the group $G_{6,2}$ is empty}.

\subsubsection{Class $P_{6,3}$}  % subsubsection 3.6.3
The algebra
$\mathcal{L}_{6,3}=L\{e_{12}-e_{14},\,e_{23}+e_{34},\,e_1,\,e_2,\,e_3,\,e_4\}$
is an extension of the algebra $\mathcal{L}_{5,3}$ by means of the vector
$e_{23}+e_{34}$, hence $P_{6,3}\subset P_{5,3}$. Substituting
\eqref{Cl:P{5,3}} for $A_i$ in \eqref{CondInv A_i:e_{23}+e_{34}(2,11a)},
we obtain the following result.
\begin{state} \label{state:P{6,3}}
The class $P_{6,3}$ of potentials that admit the group $G_{6,3}$
takes the form
\begin{equation}\label{Cl:P{6,3}}
A_i=(0,\,A,\,0,\,A),\ \ A=\text{\rm const}.
\end{equation} %\eqno(3.212)
\end{state}

\subsubsection{Class $P_{6,4}$}  % subsubsection 3.6.4
The algebra
$\mathcal{L}_{6,4}=L\{e_{12}-e_{14},\,e_{24},\,e_1,\,e_2,\,e_3,\,e_4\}$
is an extension of the algebra $\mathcal{L}_{5,1}$ by means of the vector
$e_{12}-e_{14}$, hence $C_{6,4}\subset C_{5,1}$. Substituting \eqref{Cl:P{5,1}}
for $A_i$ in \eqref{CondInv A_i:prb}--\eqref{X:prb-gl2} ($\lambda=\mu=0$),
we get the following result.
\begin{state} \label{state:P{6,4}}
The class $P_{6,4}$ of potentials that admit the group $G_{6,4}$
is defined by \eqref{Cl:P{6,3}}.
\end{state}
\begin{rem}\label{rem:P{6,4}}\rm
Statements \ref{state:P{6,3}} and \ref{state:P{6,4}} involve the potential
\eqref{Cl:P{6,3}} admits more wide symmetry group than $G_{6,3}$ and
$G_{6,4}$.
\end{rem}

\subsubsection{Class $P_{6,5}$}  % subsubsection 3.6.5
The algebra
$$
\mathcal{L}_{6,5}=L\{e_{12}-e_{14},\ e_{23}+e_{34},\ e_{13}+\lambda e_2,\
e_1,\ e_3,\ e_2-e_4\}
$$
is an extension of the algebra $\mathcal{L}_{5,5}$ ($\lambda=0$)
by means of the vector $e_{13}+\lambda e_2$, hence
$P_{6,5}\subset P_{5,5(\lambda=0)}$. Substituting \eqref{Cl:P{5,5}-0}
for $A_i$ in  \eqref{CondInv A_i-el} ($\mu=0$), we get $\lambda\Phi'=0$.
We have a result.
\begin{state} \label{state:P{6,5}}
For $\lambda=0$ the class $P_{6,5}$ of potentials that admit the group
$G_{6,5}$ is defined by \eqref{Cl:P{5,5}-0}; for $\lambda\ne0$ the same
is \eqref{Cl:P{6,3}}.
\end{state}

\subsubsection{Class $P_{6,6}$}  % subsubsection 3.6.6
The algebra
$$
\mathcal{L}_{6,6}=L\{e_{12}-e_{14},\,e_{23}+e_{34},\,e_{24},\,e_1,\,e_3,\,
e_2-e_4\}
$$
is an extension of the algebra $\mathcal{L}_{5,5}$ ($\lambda=0$)
by means of the vector $e_{24}$, hence $P_{6,6}\subset P_{5,5(\lambda=0)}$.
Substituting \eqref{Cl:P{5,5}-0} for $A_i$ in \eqref{CondInv A_i-Hyp}
($\lambda=0$), we get
\begin{equation}\label{Cl:P_{6,6}}
A_i=\left(0,\,\frac{B}{x^2+x^4},\,0,\,\frac{B}{x^2+x^4}\right),\ \
B=\text{\rm const}.
\end{equation} %\eqno(3.213)
\begin{state} \label{state:P{6,6}}
The class $P_{6,6}$ of potentials that admit the group $G_{6,6}$
is defined by \eqref{Cl:P_{6,6}}.
\end{state}

\subsubsection{Class $P_{6,7}$}  % subsubsection 3.6.7
The algebra
$$
\mathcal{L}_{6,7}=L\{e_{12}-e_{14},\,e_{23}+e_{34},\,e_{13}+\lambda e_{24},\,
e_1,\,e_3,\,e_2-e_4\}
$$
is an extension of the algebra $\mathcal{L}_{5,5}$ ($\lambda=0$)
by means of the vector $e_{13}+\lambda e_{24}$, hence
$P_{6,7}\subset P_{5,5(\lambda=0)}$. Substituting \eqref{Cl:P{5,5}-0}
for $A_i$ in \eqref{CondInv A_i:birot}--\eqref{X:birot-gl2},
we get the result.
\begin{state} \label{state:P{6,7}}
The class $P_{6,7}$ of potentials that admit the group $G_{6,7}$
is defined by \eqref{Cl:P_{6,6}}.
\end{state}
\begin{rem}\label{rem:P{6,7}}\rm
Statements \ref{state:P{6,6}} and \ref{state:P{6,7}} involve the potential
\eqref{Cl:P_{6,6}} admits more wide symmetry group than $G_{6,6}$ and
$G_{6,7}$.
\end{rem}

\subsubsection{Class $P_{6,8}$}  % subsubsection 3.6.8
The algebra
$\mathcal{L}_{6,8}=L\{e_{12},\,e_{13},\,e_{23},\,e_1,\,e_2,\,e_3\}$
is an extension of $\mathcal{L}_{3,20}$ by means of $e_1$, $e_2$,
and $e_3$, hence ${P_{6,8}\subset P_{3,20}}$. Substituting \eqref{Cl:P{3,20}}
for $A_i$ in equations $\partial_1A_i=\partial_2A_i=\partial_3A_i=0$,
we get the result.
\begin{state} \label{state:P{6,8}}
The class $P_{6,8}$ of potentials that admit the group $G_{6,8}$
takes the form $A_i=(0,\ 0,\ 0,\ \Phi(x^4))$, where $\Phi(x^4)$ is
an arbitrary function.
\end{state}

\subsubsection{Class $P_{6,9}$}  % subsubsection 3.6.9
The algebra
$\mathcal{L}_{6,9}=L\{e_{12},\ e_{14},\ e_{24},\ e_1,\ e_2,\ e_4\}$
is an extension of $\mathcal{L}_{3,21}$ by means of $e_1$, $e_2$, and $e_4$,
hence $P_{6,9}\subset P_{3,21}$. Substituting \eqref{Cl:P{3,21}}
for $A_i$ in equations $\partial_1A_i=\partial_2A_i=\partial_4A_i=0$,
we get the result.
\begin{state} \label{state:P{6,9}}
The class $P_{6,9}$ of potentials that admit the group $G_{6,9}$
takes the form $A_i=(0,\ 0,\ \Phi(x^3),\ 0)$, where $\Phi(x^3)$ is
an arbitrary function.
\end{state}
%%%%%%%%%%%%%%%%%%%%%%%%%%%%%%%%%%%%%%%%%%%%%%%%%%%%%%%%%%%%%%%%%%%%%%%%

\section{Appendix. Seven classes of Maxwell spaces that admit subgroups
of the Poincar\'{e} group.} % Section 4
\renewcommand{\theequation}{4.\arabic{equation}}
\setcounter{equation}{0}

Using the group classification of potential structures, we define more
precisely classes of Maxwell spaces that admit subgroups of the Poincar\'{e}
group \cite{Par03,Par04}. Here we describe classes of Maxwell spaces that
correspond to algebras $\mathcal L_{3,19}$, $\mathcal L_{4,16}$,
$\mathcal L_{4,17}$, $\mathcal L_{4,20}$, $\mathcal L_{5,9}$,
$\mathcal L_{6,5}$, and $\mathcal L_{6,7}$ (according to \cite{Par03}
these classes are empty).

$1^\circ.$ {\it Class $C_{3,19}$.}
For the algebra
$$
\mathcal{L}_{3,19}=L\{e_{12}-e_{14},\,e_{23}+e_{34},\,e_{13}+\lambda e_{24}\}\
\ (\lambda\ne0)
$$
we have the result.

\begin{state}\label{theor:C_{3,19}}
A Maxwell space of the class $C_{3,19}$ is defined by the tensor $F_{ij}$
such that
\begin{equation}\label{Cl:C{3,19}}
\begin{split}
&F_{12}=-\frac{\tilde x^1\Phi_2}2\left(1+\left(\tilde x^2\right)^2-
 \left(\tilde x^3\right)^2\right)-\tilde x^1\tilde x^2\tilde x^3\Phi_1+
 \frac{\Phi_3}{\tilde x^1}-\tilde x^2\Phi_5,\\
&F_{13}=\tilde x^1(\tilde x^2\Phi_1-\tilde x^3\Phi_2),\ \
 F_{14}=F_{12}+\tilde x^1\Phi_2,\\
&F_{23}=-\frac{\tilde x^1\Phi_1}2\left(1-\left(\tilde x^2\right)^2+
 \left(\tilde x^3\right)^2\right)-\tilde x^1\tilde x^2\tilde x^3\Phi_2-
 \frac{\Phi_4}{\tilde x^1}-\tilde x^3\Phi_5,\\
&F_{24}=\tilde x^1(\tilde x^2\Phi_2+\tilde x^3\Phi_1)+\Phi_5,\ \
 F_{34}=-F_{23}-\tilde x^1\Phi_1,
\end{split}
\end{equation} %\eqno(4.1)
where $\Phi_k=\Phi_k(\tilde x^4)$ $(k=1,2,5)$ are arbitrary functions,
\begin{equation}\label{C_{3,19}:Phi_3,Phi_4}
\begin{split}
&\Phi_3=\Phi_3(\tilde x^4)=\int\left(\frac1{2\lambda}\Phi_1(\tilde x^4)-
 \frac{\tilde x^4}2\Phi'_2(\tilde x^4)-
 \Phi_2(\tilde x^4)\right)d\tilde x^4,\\
&\Phi_4=\Phi_4(\tilde x^4)=\int\left(-\frac1{2\lambda}\Phi_2(\tilde x^4)+
\frac{\tilde x^4}2\Phi'_1(\tilde x^4)+
 \Phi_1(\tilde x^4)\right)d\tilde x^4,
\end{split}
\end{equation} %\eqno(4.2)
and the transformation of coordinates is defined by \eqref{Change:P_{2,11a}}.
\end{state}

\begin{example}\rm
If we replace $C_1$, $C_2$, $C_3$ by zero and $C_4$ by
$\varphi=\varphi(\tilde x^4)$ in \eqref{Solution:P_{2,11a}}--%
\eqref{Phi,Psi,Xi,Theta:P{3,19}}, we get the potential
\begin{equation}\label{Example:P{3,19}}
\begin{split}
&A_1=\varphi(\tilde x^4)\sin\frac{\ln\tilde x^1}{\lambda},\ \
 A_3=\varphi(\tilde x^4)\cos\frac{\ln\tilde x^1}{\lambda},\\
&A_2=A_4=\tilde x^2\varphi(\tilde x^4)\sin\frac{\ln\tilde x^1}{\lambda}-
 \tilde x^3\varphi(\tilde x^4)\cos\frac{\ln\tilde x^1}{\lambda}.
\end{split}
\end{equation} %\eqno(4.3)
Substituting \eqref{Example:P{3,19}}--\eqref{Change:P_{2,11a}} for $A_i$
in \eqref{F_{ij}=}, we get
\begin{align} %\label{Example:C{3,19}}
F_{12}&=-\left(\frac{2(x^1)^2\varphi'+\varphi}{x^2+x^4}+2x^2\varphi'\right)
 \sin\frac{\ln(x^2+x^4)}{\lambda}-\nonumber\\
&-\frac{2\lambda x^1x^3\varphi'+\varphi}{\lambda(x^2+x^4)}
 \cos\frac{\ln(x^2+x^4)}{\lambda},\nonumber\\
F_{13}&=2\varphi'\cdot\left(x^1\cos\frac{\ln(x^2+x^4)}{\lambda}-
x^3\sin\frac{\ln(x^2+x^4)}{\lambda}\right),\nonumber\\
F_{14}&=F_{12}+2(x^2+x^4)\varphi'\cdot\sin\frac{\ln(x^2+x^4)}{\lambda},
\nonumber\\
F_{23}&=\left(\frac{2(x^3)^2\varphi'+\varphi}{x^2+x^4}+2x^2\varphi'\right)
 \cos\frac{\ln(x^2+x^4)}{\lambda}+\nonumber\\
&+\frac{2\lambda x^1x^3\varphi'-\varphi}{\lambda(x^2+x^4)}
 \sin\frac{\ln(x^2+x^4)}{\lambda},\nonumber\\
F_{24}&=-2\varphi'\cdot\left(x^1\sin\frac{\ln(x^2+x^4)}{\lambda}+
 x^3\cos\frac{\ln(x^2+x^4)}{\lambda}\right),\nonumber
\end{align}
\begin{align}\label{Example:C{3,19}}
F_{34}=-F_{23}+2(x^2+x^4)\varphi'\cdot\cos\frac{\ln(x^2+x^4)}{\lambda}.
\end{align} %\eqno(4.4)
\end{example}
\begin{state}
If $\varphi'=\varphi'(\tilde x^4)\ne0$, then the Maxwell space defined by
the tensor \eqref{Example:C{3,19}} admits the three-di\-men\-si\-onal group
$G_S=G_{3,19}$.
\end{state}

$2^\circ.$ {\it Class $C_{4,16}$.}
For the algebra
$$
\mathcal L_{4,16}=L\{e_{12}-e_{14}+\lambda e_3,\,e_{23}+e_{34}+\lambda e_1,
\break e_{13},\,e_2-e_4\}
$$
we have the result.
\begin{state}\label{theor:C_{4,16}}
A Maxwell space of the class $C_{4,16}$ is defined by the tensor $F_{ij}$
such that
\begin{equation}\label{Cl:C{4,16}}
\begin{split}
&F_{12}=F_{14}=-\varphi\,\Phi_1(u)+\psi\,\Phi_2(u),\ \ F_{13}=\Phi_1(u),\\
&F_{23}=-F_{34}=\varphi\,\Phi_2(u)+\psi\,\Phi_1(u),\ \ F_{24}=\Phi_2(u),
\end{split}
\end{equation} %\eqno(4.5)
where $\Phi_1(u)={K}/(u^2+\lambda^2)$, $K=\text{\rm const}$, $\Phi_2(u)$ is
an arbitrary function, and
\begin{equation}\label{Change:C{4,16}}
u=x^2+x^4,\ \ \varphi=\frac{\lambda x^1+ux^3}{u^2+\lambda^2},\ \
\psi=\frac{\lambda x^3-x^1u}{u^2+\lambda^2}.
\end{equation} %\eqno(4.6)
\end{state}
\begin{example}\rm
Substituting $0$ for $\Phi_2(u)$ in \eqref{Cl:C{4,16}}--\eqref{Change:C{4,16}},
we get
\begin{equation}\label{Example:C{4,16}}
\begin{split}
&F_{12}=F_{14}= -\frac{K(\lambda x^1+(x^2+x^4)x^3)}
 {\left((x^2+x^4)^2+\lambda^2\right)^2},\ \
 F_{13}=\frac{K}{(x^2+x^4)^2+\lambda^2},\\
&F_{23}=-F_{34}=\frac{K(\lambda x^3-x^1(x^2+x^4))}
 {\left((x^2+x^4)^2+\lambda^2\right)^2},\ \ F_{24}=0.
\end{split}
\end{equation} %\eqno(4.7)
\end{example}
\begin{state}
If $K\ne0$, then the Maxwell space defined by the tensor
\eqref{Example:C{4,16}} admits the four-di\-men\-si\-onal group
$G_S=G_{4,16}$.
\end{state}

$3^\circ.$ {\it Class $C_{4,17}$.}
For the algebra
$$
\mathcal{L}_{4,17}=L\{e_{12}-e_{14},\,e_{23}+e_{34},\,e_{13}+\lambda e_{24},\,
e_2-e_4\}
$$
we have the result.
\begin{state}\label{theor:C_{4,17}}
A Maxwell space of the class $C_{4,17}$ is defined by the tensor $F_{ij}$
such that
\begin{align}\label{Cl:C{4,17}}
&F_{12}=F_{14}=\frac{1}{x^2+x^4}\left(A\cos\frac{\ln(x^2+x^4)}{\lambda}+
 B\sin\frac{\ln(x^2+x^4)}{\lambda}+Cx^1\right),\nonumber\\
&F_{23}=-F_{34}=\frac{1}{x^2+x^4}\left(B\cos\frac{\ln(x^2+x^4)}{\lambda}-
 A\sin\frac{\ln(x^2+x^4)}{\lambda}-Cx^3\right),\nonumber\\
&F_{13}=0,\ \ F_{24}=C\ \ (A,\,B,\,C=\text{\rm const}).
\end{align} %\eqno(4.8)
\end{state}
\begin{state}
If 1) $C\ne0$ and 2) $A\ne0$ (or $B\ne0$), then the Maxwell space
defined by the tensor \eqref{Cl:C{4,17}} admits the four-di\-men\-si\-onal
group $G_S=G_{4,17}$.
\end{state}

$4^\circ.$ {\it Class $C_{4,20}$.}
For the algebra
$\mathcal{L}_{4,20}=L\{e_{12}-e_{14},\,e_{23}+e_{34},\,e_{13},\,e_{24}\}$
we have the result.
\begin{state}\label{theor:C_{4,20}}
A Maxwell space of the class $C_{4,20}$ is defined by the tensor $F_{ij}$
such that
\begin{equation}\label{Cl:C{4,20}}
\begin{split}
&F_{12}=F_{14}=\frac{x^1\Phi}{x^2+x^4},\ \ F_{13}=0,\ \
 F_{23}=-F_{34}=-\frac{x^3\Phi}{x^2+x^4},\\
&F_{24}=\Phi\ \ (\Phi=\Phi(\tilde x^4)=\Phi((x^1)^2+(x^2)^2+(x^3)^2-(x^4)^2)).
\end{split}
\end{equation} %\eqno(4.9)
\end{state}
\begin{state}
If $\Phi'(\tilde x^4)\ne0$, then the Maxwell space defined by the tensor
\eqref{Cl:C{4,20}} admits the four-di\-men\-si\-onal group $G_S=G_{4,20}$.
\end{state}
\begin{rem}
In \cite{Par04} the class $C_{4,20}$ is not empty, but this is more narrow
than above.
\end{rem}

$5^\circ.$ {\it Class $C_{5,9}$.}
For the algebra
$$
\mathcal{L}_{5,9}=L\{e_{12}-e_{14},\ e_{23}+e_{34},\ e_{13},\
e_{24},\ e_2-e_4\}
$$
we have the result.
\begin{state}\label{theor:C_{5,9}}
A Maxwell space of the class $C_{5,9}$ is defined by the tensor $F_{ij}$
such that
\begin{equation}\label{Cl:C{5,9}}
\begin{split}
&F_{12}=F_{14}=\frac{Cx^1}{x^2+x^4},\ \ F_{13}=0,\ \ F_{24}=C,\\
&F_{23}=-F_{34}=-\frac{Cx^3}{x^2+x^4},\ \ (C=\text{\rm const}).
\end{split}
\end{equation} %\eqno(4.10)
\end{state}
\begin{state}
If $C\ne0$, then the Maxwell space defined by the tensor \eqref{Cl:C{5,9}}
admits the five-di\-men\-si\-onal group $G_S=G_{5,9}$.
\end{state}

$6^\circ.$ {\it Class $C_{6,5}$.}
For the algebra
$$
\mathcal{L}_{6,5}=L\{e_{12}-e_{14},\ e_{23}+e_{34},\ e_{13}+\lambda e_2,\
e_1,\ e_3,\ e_2-e_4\}
$$
we have the result.
\begin{state}\label{theor:C_{6,5}}
For ${\lambda\ne0}$ a Maxwell space of the class $C_{6,5}$ is defined by
the tensor $F_{ij}$ such that
\begin{equation}\label{Cl:C{6,5}}
\begin{split}
&F_{12}=F_{14}=C_1\sin\frac{x^2+x^4}{\lambda}+C_2\cos\frac{x^2+x^4}{\lambda},\\
&F_{23}=-F_{34}=C_1\cos\frac{x^2+x^4}{\lambda}-C_2\sin\frac{x^2+x^4}{\lambda},\\
&F_{13}=F_{24}=0\ \ (C_1,\,C_2=\text{\rm const});
\end{split}
\end{equation} %\eqno(4.11)
if $\lambda=0$, then $F_{ij}=0$.
\end{state}
\begin{state}
If $C_1\ne0$ or $C_2\ne0$, then the Maxwell space defined by the tensor
\eqref{Cl:C{6,5}} admits the six-di\-men\-si\-onal group $G_S=G_{6,5}$.
\end{state}

$7^\circ.$ {\it Class $C_{6,7}$.}
For the algebra
$$
\mathcal{L}_{6,7}=L\{e_{12}-e_{14},\,e_{23}+e_{34},\,e_{13}+\lambda e_{24},\,
e_1,\,e_3,\,e_2-e_4\}
$$
we have the result.
\begin{state}\label{theor:C_{6,7}}
A Maxwell space of the class $C_{6,7}$ is defined by the tensor $F_{ij}$
such that
\begin{equation}\label{Cl:C{6,7}}
F_{12}=F_{14}=\Phi,\ \ F_{13}=F_{24}=0,\ \ F_{23}=-F_{34}=\Psi,
\end{equation} %\eqno(4.12)
where
\begin{equation}\label{Cl:C{6,7}-Phi,Psi}
\begin{split}
&\Phi=\frac1{x^2+x^4}\left(a_1\cos\frac{\ln(x^2+x^4)}{\lambda}-
 a_2\sin\frac{\ln(x^2+x^4)}{\lambda}\right),\\
&\Psi=\frac1{x^2+x^4}\left(a_1\sin\frac{\ln(x^2+x^4)}{\lambda}+
 a_2\cos\frac{\ln(x^2+x^4)}{\lambda}\right)
\end{split}
\end{equation} %\eqno(4.13)
($a_1,\,a_2=\text{\rm const})$.
\end{state}
\begin{state}
If $a_1\ne0$ or $a_2\ne0$, then the Maxwell space defined by the tensor
\eqref{Cl:C{6,7}}--\eqref{Cl:C{6,7}-Phi,Psi} admits the six-di\-men\-si\-onal
group ${G_S=G_{6,5}}$.
\end{state}

%\end{document}

%%%%%%%%%%%%%%%%%%%%%%%%%%%%%%%%%%%%%%%%%%%%%%%%%%%%%%%%%%%%%%%%%%%%%%%%

\end{document}